\newcommand{\R}{I\!\!R}
\newcommand{\Z}{\mathbb Z}

\newcommand{\M}{\mathcal M}
\newcommand{\g}{\mathfrak g}
\newcommand{\Pc}{\mathcal P}
\newcommand{\Qc}{\mathcal Q}
\newcommand{\J}{\mathfrak J}
\newcommand{\Jb}{\mathbb J}
\newcommand{\Rc}{\mathcal R}
\newcommand{\Sc}{\mathcal S}

\newcommand{\SgD}{\mathfrak S_\gamma^{\mathcal D}}
\newcommand{\KgD}{\mathcal K_{\gamma,\Pc}^{\mathcal D}}
\newcommand{\KgDPQ}{\mathcal K_{\gamma,\Pc,\Qc}^{\mathcal D}}
\newcommand{\St}{\mathfrak S}
\newcommand{\Kt}{\mathcal K}

\newcommand{\sgn}{\mathrm{sgn}}
\newcommand{\dgn}{\mathrm{dgn}}
\newcommand{\rk}{\mathrm{rk}}
\newcommand{\mul}{\mathrm{mul}}
\newcommand{\ifoc}{\mathrm i_{\mathrm{foc}}}
\newcommand{\imaslov}{\mathrm i_{\mathrm{maslov}}}

\newcommand{\Bsym}{\mathrm B_{\mathrm{sym}}}
\newcommand{\Iota}[2]{\mathcal D_{#1,#2}}
\newcommand{\Lge}{\Lambda_{\ge1}}

\documentclass[oneside,draft,letterpaper,10pt,titlepage]{amsart}

\usepackage{amsmath}               
\usepackage{amsthm}                
\usepackage{times}
\usepackage{amscd}

\numberwithin{equation}{section}

\title[The semi-Riemannian Morse Index Theorem]%
{The Morse Index Theorem in  semi-Riemannian Geometry}
\author[P.\ Piccione]{Paolo Piccione}
\author[D.\ Tausk]{Daniel V.\ Tausk}
\address{Departamento de Matem\'atica,\hfill\break\indent
Instituto de Matem\'atica e Estat\'\i stica\hfill\break\indent  Universidade de S\~ao
Paulo, \hfill\break\indent
Caixa Postal 66281, CEP 05315--970, SP\hfill\break\indent Brazil}
\email{piccione@ime.usp.br, tausk@ime.usp.br}
\urladdr{http://www.ime.usp.br/\~{}piccione, http://www.ime.usp.br/\~{}tausk}
\thanks{The first author is partially sponsored by CNPq (Processo n.\ 301410/95), the
second author is sponsored by FAPESP (Processo n.\ 98/12530-2)}

\date{November 2000}

\begin{document}


\theoremstyle{plain}\newtheorem{teo}{Theorem}[section]
\theoremstyle{plain}\newtheorem{prop}[teo]{Proposition}
\theoremstyle{plain}\newtheorem{lem}[teo]{Lemma} 
\theoremstyle{plain}\newtheorem{cor}[teo]{Corollary} 
\theoremstyle{definition}\newtheorem{defin}[teo]{Definition} 
\theoremstyle{remark}\newtheorem{rem}[teo]{Remark} 
\theoremstyle{plain} \newtheorem{assum}[teo]{Assumption}
\theoremstyle{definition}\newtheorem{example}[teo]{Example}


\begin{abstract}
We prove a semi-Riemannian version of the celebrated Morse Index
Theorem for geodesics in semi-Riemannian manifolds; we consider
the general case of both endpoints variable on two submanifolds.
The key role of the theory is played by the notion of the {\em Maslov index\/}
of a semi-Riemannian geodesic, which is a homological invariant
and it substitutes the notion of geometric index in Riemannian
geometry. Under generic circumstances, the Maslov index of a geodesic
is computed as a sort of {\em algebraic count\/} of the
conjugate points along the geodesic. For non positive definite metrics
the index of the index form is always infinite; in this paper we prove
that the space of all variations of a given geodesic has a {\em natural\/} 
splitting into two infinite dimensional subspaces, and the Maslov index 
is given by the difference of the index and the coindex of the restriction of
the index form to these subspaces. In the case of variable endpoints,
two suitable correction terms, defined in terms of the endmanifolds, 
are added to the equality. Using appropriate change of variables,
the theory is entirely extended to the more general case of
{\em symplectic differential systems}, that can be obtained as linearizations
of the Hamilton  equations. The main results proven in this paper were announced
in \cite{PT3}.
\end{abstract}

\maketitle

\begin{section}{Introduction}\label{sec:intro}
Let $(\mathcal M,\g)$ be a Riemannian manifold; the classical
 Morse Index Theorem  states
that the number of conjugate points along a geodesic $\gamma:[a,b]\to\mathcal M$ counted
with multiplicities (the geometric index of $\gamma$) is equal to the  index of the
second variation of the Riemannian action  functional $E(z)\!=\!\frac12\int_a^b \g(\dot
z,\dot z)\,\mathrm dt$ at the critical point $\gamma$. Such second variation is called
the {\em index form}, and it will be  denoted by
$I_\gamma$. The theorem has later been extended in several directions (see~\cite{BE,
BEE, dC, Duis, Ed, EK, Kal, M, ON, Sm} for   versions of this theorem in different
contexts). In Lorentzian geometry, the theorem holds in the case of causal (i.e., 
nonspacelike) geodesics, provided that one considers the restriction of $I_\gamma$ to
the space of variations that are everywhere orthogonal to the geodesic. However, when
one considers the case of spacelike Lorentzian geodesics or geodesics in semi-Riemannian
manifolds with metric of arbitrary  index, there is no hope to extend the original
formulation of the theorem, due mainly to the following phenomena:
\begin{itemize}
\item the set of conjugate points along a geodesic may fail to be
discrete (see \cite{Hel1, fechado});
\item the index of $I_\gamma$ is always infinite, even when restricted to the
space of variations orthogonal to $\gamma$ (see Proposition~\ref{thm:indinfinito}).
\end{itemize}
The case of spacelike Lorentzian geodesics has been studied
in \cite{GMPT}, where the authors consider a {\em stationary\/}
metric $\g$, i.e., a metric admitting a timelike Killing vector field
$Y$.  The Killing field $Y$ gives a conservation law for geodesics $\gamma$:
$\g(\dot\gamma,Y)=\text{constant}$;
the main result of the paper is that, if one restricts the index
form to the space of variational vector fields along $\gamma$ corresponding
to variations of $\gamma$ by curves that satisfy
such conservation law, then the index of this restriction is finite, and it is equal to
a homological invariant of the geodesic called the {\em Maslov index}.
The notion of Maslov index associated to curves in a Lagrangian submanifold of
$\R^{2n}$  appeared originally in the Russian
literature (see for instance \cite{Ar} and the references therein).
Some interesting applications in Variational Calculus 
of the   Maslov index were shown by Duistermaat in \cite{Duis},
where it is proven an index theorem for solutions of convex Hamiltonian
systems. 
An index theorem for solutions of non convex Hamiltonian systems is 
proven in \cite{PT2}; the result of \cite{PT2} is a weak form of the index theorem
proven in this paper in a sense clarified below.

There is nowadays quite an extensive literature concerning applications
of the Maslov index to the theory of Hamiltonian systems (see for instance
\cite{CZ, Long, SZ});
in the context of semi-Riemannian geodesics the Maslov index 
 was introduced by Helfer in \cite{Hel1}. Under a suitable nondegeneracy assumption, that holds {\em generically},
one proves that each conjugate point along a semi-Riemannian geodesic
is isolated, and that the Maslov index of the geodesic is given by the 
sum of the {\em signatures\/} of the conjugate points (see
Definition~\ref{thm:defsignature}). The Maslov index is defined in general
as the intersection number of a curve $\ell$ in the Lagrangian Grassmannian $\Lambda$
of a symplectic space with the codimension one, transversally oriented
subvariety of $\Lambda$, consisting of those Lagrangians that are not
transverse to a fixed one. The curve $\ell$ is obtained from the flow of
the Jacobi equation along $\gamma$.

The main purpose of this paper is to determine the relations
between the Maslov index of a semi-Riemannian geodesic $\gamma$ and
the index form $I_\gamma$, obtaining a general version of
the Morse index theorem in semi-Riemannian geometry. More precisely, 
generalizing the ideas
in \cite{GMPT, PT2}, we prove that the choice of a {\em maximal negative distribution\/}
along
$\gamma$ determines  a {\em natural\/}  splitting of the space of
all variations  of $\gamma$ into two $I_\gamma$-orthogonal infinite dimensional 
subspaces ${\mathcal K_\gamma}$, $\mathcal S_\gamma$ such that the Maslov index  is given
by the difference of the index of
$I_\gamma\vert_{\mathcal K_\gamma}$ and the coindex of $I_\gamma\vert_{\mathcal
S_\gamma}$ (i.e., the index of $-I_\gamma\vert_{\mathcal S_\gamma}$).
 This kind of result  aims to a generalized Morse theory
for strongly indefinite functionals on Hilbert manifolds (see \cite{Abb}).

A different index theory for semi-Riemannian geodesics   is 
presented in \cite{Hel1}, where, under a suitable nondegeneracy assumption, the author
proves an equality between the Maslov index  and the {\em spectral index\/} of the
geodesic, which is an integer number defined in terms of the spectral properties
of the Jacobi differential operator.
Also, in \cite{Hel1} there is an attempt to relate the spectral index
with the difference between the index and the coindex of suitable restrictions
of $I_\gamma$. However, the construction discussed by Helfer has no geometrical
interpretation, and, as a matter of facts, it is not hard to prove that, by minor
modifications of this construction, one can produce {\em any\/} integer number as a
difference of the index and the coindex of restrictions of $I_\gamma$.
A further discussion of Helfer's results can be found in references \cite{MPT, PT4}.

In order to motivate the main result of this paper, we can consider the following simple
but instructive example. Consider the case of a product semi-Riemannian manifold
$\mathcal M=\mathcal M_1\times \mathcal M_2$, endowed with the
metric $\g=\g_1\oplus(-\g_2)$, where $\g_1$, $\g_2$ are Riemannian metrics on
$\mathcal M_1$ and $\mathcal M_2$ respectively. If $\gamma=(\gamma_1,\gamma_2)$ is
a geodesic in $\mathcal M$, the set of conjugate points along $\gamma$ is given by
the union of the set of conjugate points along $\gamma_1$ and the set of conjugate
points along $\gamma_2$. Using the Riemannian Morse Index Theorem it is easily seen that
the index of the restriction of $I_\gamma$ to the space $\mathcal K_\gamma$ of
variational vector fields along $\gamma_1$ equals the number of conjugate points along
$\gamma_1$, while the coindex of the restriction of $I_\gamma$ (i.e., the index of
$-I_\gamma$) to the space $\mathcal S_\gamma$ of variational vector fields along
$\gamma_2$ equals the number of conjugate points along $\gamma_2$. In this case, the
Maslov index of $\gamma$ equals the difference between the geometric indexes of
$\gamma_1$ and
$\gamma_2$.

The idea of the construction of the spaces $\mathcal K_\gamma$ and $\mathcal S_\gamma$
for the general case is the following.
One considers a maximal distribution $\mathcal D$ of subspaces along the geodesic
$\gamma$ on which the metric is negative definite; in the above example, $\mathcal D$
would be given by $T\mathcal M_2$.
The space $\mathcal S_\gamma$ is defined as the space of variational vector fields along
$\gamma$ taking values in $\mathcal D$. The space $\mathcal K_\gamma$ is defined as the
space of variational vector fields along $\gamma$ that are {\em Jacobi in the directions
of $\mathcal D$}, that is, vector fields whose image by the Jacobi differential
operator is orthogonal to the distribution $\mathcal D$.  One proves
that the restrictions of $I_\gamma$ to $\mathcal S_\gamma$ and $\mathcal K_\gamma$ are
represented by a compact perturbation of a negative and a positive
isomorphism, respectively, and therefore $n_+\big(I_\gamma\vert_{\mathcal S_\gamma}\big)$
and $n_-\big(I_\gamma\vert_{\mathcal K_\gamma}\big)$ are finite natural numbers.
Here, by $n_-$ and $n_+$ we mean respectively the index and the coindex of
a symmetric bilinear form. 

The spaces $\mathcal K_\gamma$ and $\mathcal S_\gamma$ are naturally associated to the
quadruple
$(\mathcal M,\g,\gamma,\mathcal D)$ in the following categorical sense. If $F:(\mathcal
M, \g) \to (\widetilde{\mathcal M} ,\tilde{\g})$ is an isometry sending 
$\gamma$ to $\tilde\gamma$ and $\mathcal
D$ onto $\tilde{\mathcal D}$, then $F$ also sends the   spaces $\mathcal
K_\gamma,\mathcal S_\gamma$ corresponding to $(\mathcal M,\g,\gamma,\mathcal D)$ to the
spaces
$\tilde{\mathcal K_\gamma}$ and $\tilde{\mathcal S_\gamma}$ corresponding to 
$(\widetilde{\mathcal M},\tilde{\g},\tilde\gamma,\tilde{\mathcal D})$.

Let us now give a brief description of the technique used to prove our
main result.

The computation of $n_+\big(I_\gamma\vert_{\mathcal S_\gamma}\big)$ is done by 
proving that $-I_\gamma\vert_{\mathcal S_\gamma}$ corresponds to the index 
form of a  {\em positive definite symplectic
system\/} (Subsections~\ref{sub:symplectic} and \ref{sub:reduced}); in this case the
classical Morse Index Theorem applies.

The computation of the index $n_-\big(I_\gamma\vert_{\mathcal K_\gamma}\big)$
is done by considering the evolution of the function
$i(t)=n_-\big(I_\gamma(t)\vert_{\mathcal K_\gamma(t)}\big)$, 
where $I_\gamma(t)$ is the index form of the restriction $\gamma\vert_{[a,t]}$
and $\mathcal K_\gamma(t)$ is the corresponding restricted version of $\mathcal
K_\gamma$. By a perturbation argument, one can assume that there is only a finite number
of conjugate points along $\gamma$, in which case $i$ is piecewise constant
(although not necessarily monotonic). The jumps of $i$ occur at those
instants $t$ for which $\gamma(t)$ is  conjugate and also when $\mathcal K_\gamma(t)\cap
\mathcal S_\gamma(t)\ne\{0\}$; here, by $\mathcal S_\gamma(t)$ we mean the restricted
version of the space $\mathcal S_\gamma$.  
In studying the evolution of the function $i$, a technical problem arises
due to the fact that the family $\mathcal K_\gamma(t)$ does not vary smoothly with
respect to $t$; indeed, the family may have singularities at those
instants $t$ when $\mathcal K_\gamma(t)\cap \mathcal S_\gamma(t)\ne\{0\}$. 
In order to overcome this problem,   we introduce an auxiliary extension
$I^\#_\gamma(t)$ of the index form and an auxiliary extension $\mathcal K^\#_\gamma(t)$
of $\mathcal K_\gamma(t)$ such that:
\begin{itemize}
\item $\mathcal K^\#_\gamma(t)$ varies smoothly with $t$;
\item for $t\ne t_0$, the indexes of $I^\#_\gamma(t)\vert_{\mathcal K^\#_\gamma(t)}$ and
of $I_\gamma(t)\vert_{\mathcal K_\gamma(t)}$ are  easily related;
\item $I^\#_\gamma(t_0)$ is nondegenerate on $\mathcal K^\#_\gamma(t)$, therefore its
index is constant around $t=t_0$.
\end{itemize}
Using a symplectic geometry result (Lemma~\ref{thm:tech}), we conclude that
the jump of $i$ at each conjugate point coincides with its contribution to
the Maslov index of the geodesic. It is a surprising fact that virtually
all the previous versions of the Morse Index Theorem can be deduced
as a simple consequence of this Lemma. As to the jumps of $i$ corresponding to
those $t$'s for which $\mathcal K_\gamma(t)\cap\mathcal S_\gamma(t)\ne\{0\}$, we employ
a functional analytical technique which says essentially that the
jump of the index of a $C^1$ curve of symmetric bilinear forms passing
through a degenerate instant is given by the signature of the derivative
restricted to the kernel. 

For the sake of completeness, in the paper
we will consider the more general case that the initial
endpoint of the geodesic is left free to move in a nondegenerate submanifold
$\mathcal P$ of $\mathcal M$, and the notion of conjugate
point is replaced by that of $\mathcal P$-focal point. For this case, the theory
is perfectly analogous to the case of a fixed initial point, with the only
exception that the initial value of the function $i$ is in general non zero,
but it is given by the index of the restriction of the metric $\g$ to
$T_{\gamma(a)}\mathcal P$. This is an entirely new phenomenon, that can only
occur in manifolds with a nonpositive definite metric.

The index theorem in the even more general case of a geodesic with final 
endpoint varying in a submanifold $\mathcal Q$ of $\M$ is then easily
obtained by a simple observation, that appears already in \cite{PT}.
What is interesting to remark here is that this observation led the authors
to the idea of considering the auxiliary extension of the index form $I^\#_\gamma$
that was mentioned above. Namely, $I^\#_\gamma$ can be thought of as the index 
form corresponding to the geodesic $\gamma$ when the final endpoint varies in a 
fictitious submanifold. 

We outline briefly the structure of the paper.

In Section~\ref{sec:semiriemannian} we give the basic definitions
concerning focal points and the index form and in Section~\ref{sec:maslov}
we define the Maslov index. 
In Subsection~\ref{sub:Maslovcurve} we study curves in the Lagrangian Grassmannian
of a symplectic space and give a few technical lemmas in symplectic
geometry. 
In Subsection~\ref{sub:Maslovgeo} we define the Maslov index of a semi-Riemannian
geodesic. 

In Section~\ref{sec:abstract} we give some abstract functional analytical
results concerning the variation of the index of a curve of symmetric
bilinear forms on a Hilbert space.

Our main results are stated in Section~\ref{sec:index}; the proofs
are spread throughout the following subsections.
In Subsection~\ref{sub:reduction}, by means of a parallel trivialization of
the tangent bundle along the geodesic, we reduce   the problem
to the theory of Morse--Sturm systems in $\R^n$.
In Subsection~\ref{sub:symplectic} we introduce the class of symplectic
differential systems needed in the computation of the coindex $n_+\big(I_\gamma\vert_{\mathcal
S_\gamma}\big)$; the class of symplectic differential systems  extends
naturally the class of Morse--Sturm systems.
 In Subsection~\ref{sub:reduced} we introduce the {\em reduced
symplectic system}, which is naturally associated to the choice of
the maximal negative distribution $\mathcal D$.
In Subsection~\ref{sub:sustenido} we define the auxiliary extension
$I_\gamma^\#$ and we discuss its properties.
The index function $i(t)$ is introduced in Subsection~\ref{sub:indexfunction}
and in Subsection~\ref{sub:proof} we conclude the proof of our main theorem.

In Section~\ref{sec:indexthsympl}, using the fact that every symplectic system
is isomorphic to a Morse--Sturm system, we extend the theory to this context
and we obtain an Index Theorem for solutions of Hamiltonian systems.
A preliminary version of this theorem appears in \cite{PT2}, where the result
is proven under the restrictive assumption that $I_\gamma$ is negative
definite~in~$\mathcal S_{\gamma}$.
\end{section}
\begin{section}{Semi-Riemannian Geodesics}\label{sec:semiriemannian}
In this section we give the basic definitions concerning the geometry
of semi-Riemann\-ian manifolds and their geodesics.

We start with some general definitions 
concerning symmetric bilinear forms for later use.
 Let $V$ be any real vector space and $B:V\times V\to\R$ a symmetric
bilinear form; given a subspace $W\subset V$, we will denote with $B\vert_W$ the
restriction of $B$ to $W\times W$. The {\em negative type number\/} (or {\em index})
$n_-(B)$ of
$B$ is the possibly infinite number defined by
\begin{equation}\label{eq:def-}
n_-(B)=\sup\Big\{{\rm dim}(W):W\ \text{subspace of}\ V \ \text{such that}\ B\vert_W\ 
\text{is negative definite}\Big\}.
\end{equation}
The {\em positive type number\/} $n_+(B)$ (or {\em coindex}) is given by
$n_+(B)=n_-(-B)$; if at least one of these two numbers is finite,
the {\em signature\/} $\sgn(B)$ is defined by:
\[\sgn(B)=n_+(B)-n_-(B).\]
The {\em kernel\/} of $B$, ${\mathrm Ker}(B)$, is the set  of vectors $v\in V$
such that $B(v,w)=0$ for all $w\in V$; the {\em degeneracy\/} $\dgn(B)$
of $B$ is the (possibly infinite) dimension of ${\rm Ker}(B)$. 
If $V$ is finite dimensional, then the numbers $n_+(B)$, $n_-(B)$ and
$\dgn(B)$ are respectively  the number of $1$'s, $-1$'s and $0$'s
in the canonical form of $B$ as given by the Sylvester's Inertia Theorem.
In this case, $n_+(B)+n_-(B)$ is equal to the codimension of ${\rm Ker}(B)$,
and it is also called the {\em rank} of $B$, $\rk(B)$.

Let $(\M,\g)$ be an $n$-dimensional semi-Riemannian manifold, with $\g$ a metric
tensor of (constant) index~$k$:
\begin{equation}\label{eq:indexg}
n_-(\g)=k.
\end{equation}
Let $\nabla$ denote the Levi--Civita connection of $\g$ and let $\Rc$
be the corresponding curvature tensor, chosen with the following
sign convention:
\[\Rc(X,Y)=\nabla_X\nabla_Y-\nabla_Y\nabla_X-\nabla_{[X,Y]}.\]
Let $\Pc\subset\M$ be a smooth submanifold and $\gamma:[a,b]\to\M$ be
a geodesic with $\gamma(a)\in\Pc$ and $\dot\gamma(a)\in T_{\gamma(a)}\Pc^\perp$,
where $\perp$ denotes the orthogonal complement with respect to $\g$.

We assume that $\Pc$ is {\em nondegenerate\/} at $\gamma(a)$, i.e., 
that the restriction of $\g$ to $T_{\gamma(a)}\Pc$ is nondegenerate.
For $p\in\Pc$ and $n\in T_p\Pc^\perp$, the {\em second fundamental form\/}
$\Sc^\Pc_{n}$ is the symmetric bilinear form on $T_p\Pc$ defined by:
\[\Sc^\Pc_n(v_1,v_2)=\g(\nabla_{v_1}V_2,n),\]
where $V_2$ is any smooth vector field in $\Pc$ with $V_2(p)=v_2$.
Since $\Pc$ is nondegenerate at $\gamma(a)$, then $\Sc^\Pc_{\dot\gamma(a)}$ can be
thought of as a $\g$-symmetric linear endomorphism of $T_{\gamma(a)}\Pc$.

A {\em Jacobi field\/} along $\gamma$ is a smooth vector field $J$
along $\gamma$ satisfying the second order linear differential equation:
\[J''=\Rc(\dot\gamma,J)\,\dot\gamma,\]
where the prime means covariant derivative along $\gamma$. A {\em $\Pc$-Jacobi
field\/} is a Jacobi field satisfying the initial conditions:
\begin{equation}\label{eq:IC}
J(a)\in T_{\gamma(a)}\Pc,\quad\text{and}\quad J'(a)+\Sc^\Pc_{\dot\gamma(a)}(J(a))
\in T_{\gamma(a)}\Pc^\perp.
\end{equation}
We denote by $\J$ the vector space of all $\Pc$-Jacobi fields
along $\gamma$:
\begin{equation}\label{eq:J}
\J=\Big\{J: J\ \text{is $\Pc$-Jacobi along $\gamma$}\Big\};
\end{equation}
$\J$ is an $n$-dimensional vector space; for all $t\in[a,b]$, we set:
\begin{equation}\label{eq:Jt}
\J[t]=\Big\{J(t):J\in\J\Big\}\subset T_{\gamma(t)}\M.
\end{equation}
A point $\gamma(t)$, with $t\in\left]a,b\right]$, is said to be {\em $\Pc$-focal\/}
if there exists a non zero $J\in\J$ such that $J(t)=0$. 
We have that
$\gamma(t)$ is $\Pc$-focal if and only if $\J[t]\ne T_{\gamma(t)}\M$.
The {\em multiplicity\/} $\mul(t)$ of the $\Pc$-focal point $\gamma(t)$ is the
dimension of the space of those $J\in\J$ such that $J(t)=0$; the
multiplicity of $\gamma(t)$ coincides with the codimension of
$\J[t]$ in $T_{\gamma(t)}\M$. 

For non positive definite metrics, we have a more appropriate
notion of ``size'' for a $\Pc$-focal point:
\begin{defin}\label{thm:defsignature}
The {\em signature\/} $\sgn(t)$ of a $\Pc$-focal point $\gamma(t)$ is 
the signature of the restriction of $\g$ to $\J[t]^\perp$:
\[\sgn(t)=\sgn\big(\g\vert_{\J[t]^\perp}\big).\]
The $\Pc$-focal point $\gamma(t)$ is said to be {\em nondegenerate\/}
if such restriction is nondegenerate. If there are only a finite
number of $\Pc$-focal points along $\gamma$, then we define the
{\em focal index\/} $\ifoc(\gamma)$ of $\gamma$ as the sum of the signatures
of all the $\Pc$-focal points along $\gamma$:
\[\ifoc(\gamma)=\sum_{t\in\left]a,b\right]}\sgn(t).\]
\end{defin}
For instance, if $(\M,\g)$ is Riemannian ($k=0$), or
if $(\M,\g)$ is Lorentzian ($k=1$) and $\gamma$ is causal, i.e.,
$\g(\dot\gamma,\dot\gamma)\le0$, then all the $\Pc$-focal points are
nondegenerate, and their signatures coincide  with their multiplicity.
Namely, in this case $\g$ is positive definite in $\J[t]^\perp$.  

In the Riemannian or in the causal Lorentzian case it is well
known that the set of $\Pc$-focal points along a geodesic
is discrete;
in the general semi-Riemannian case,  focal points {\em may\/}
indeed accumulate (see \cite{Hel1, fechado}) even in the case that $\Pc$ is
a point. We have the
following result concerning the distribution of $\Pc$-focal points:
\begin{prop}\label{thm:isolPfocalpts}
There are no $\Pc$-focal points $\gamma(t)$ for $t$ near $a$.
Nondegenerate $\Pc$-focal points are isolated. Moreover, if $(\M,\g)$
is real analytic, then the set of $\Pc$-focal points along $\gamma$ is
finite.
\end{prop}
\begin{proof}
See for instance \cite[Proposition~2.5.1, Remark~2.5.3]{MPT}.
\end{proof}
We consider the following symmetric bilinear form:
\begin{equation}\label{eq:defIP}
I_\gamma^\Pc(v,w)=\int_a^b\Big[\g(v',w')+\g\big(\Rc(\dot\gamma,v)\,\dot\gamma,w\big)\Big]\;\mathrm
dt-\mathcal S^\Pc_{\dot\gamma(a)}\big(v(a),w(a)\big),
\end{equation}
defined on the space $\mathcal H^P_\gamma$ of all vector fields
$v$ along $\gamma$ of $H^1$-Sobolev regularity\footnote{%
this means that $v:[a,b]\to T\M$ is absolutely continuous, and the
covariant derivative $v'$ is square-integrable with respect to some
positive definite inner product along $\gamma$.}
with $v(a)\in T_{\gamma(a)}\Pc$
and $v(b)=0$. The space $\mathcal H_\gamma^\Pc$ has the topology
of a Hilbertable space, and $I_\gamma^\Pc$ is continuous in this topology.
The set $\Omega_{\Pc,\gamma(b)}$ of all curves of $H^1$-regularity in 
$\M$ joining $\Pc$ and
$\gamma(b)$ can be given the structure of an infinite dimensional Hilbert manifold,
and the semi-Riemannian {\em action functional\/} $f(z)=\frac12\int_a^b\g(\dot z,\dot z)
\;\mathrm dt$ is smooth on $\Omega_{\Pc,\gamma(b)}$. The geodesic $\gamma$
is a critical point of $f$ in $\Omega_{\Pc,\gamma(b)}$,  $\mathcal H_\gamma^\Pc$
is the tangent space $T_\gamma\Omega_{\Pc,\gamma(b)}$ and the symmetric bilinear form
$I^\Pc_\gamma$ is the {\em Hessian\/} of $f$ at $\gamma$.

We now consider another smooth submanifold $\Qc\subset\M$ with $\gamma(b)\in\Qc$
and $\dot\gamma(b)\in T_{\gamma(b)}\Qc^\perp$. In this situation, $\gamma$
is also a critical point for the action functional $f$ defined in the
Hilbert manifold $\Omega_{\Pc,\Qc}$ of all $H^1$-curves joining $\Pc$ and $\Qc$. 
The tangent space $T_\gamma \Omega_{\Pc,\Qc}$ will be denoted by
$\mathcal H_{\Pc,\Qc}$, and it consists of all $H^1$-vector fields $v$ along
$\gamma$ with $v(a)\in T_{\gamma(a)}\Pc$ and $v(b)\in T_{\gamma(b)}\Qc$.
The Hessian of $f$ at $\gamma$ in the space $\Omega_{\Pc,\Qc}$ is given by
the following bounded symmetric bilinear form in $\mathcal H_{\Pc,\Qc}$:
\begin{equation}\label{eq:defIPQ}
\begin{split}
I_\gamma^{\Pc,\Qc}(v,w)=&\int_a^b\Big[\g(v',w')+\g\big(\Rc(\dot\gamma,v)\,\dot\gamma,w\big)\Big]\;\mathrm
dt \\ &+\mathcal S^\Qc_{\dot\gamma(b)}\big(v(b),w(b)\big)-\mathcal
S^\Pc_{\dot\gamma(a)}\big(v(a),w(a)\big).
\end{split}
\end{equation}
If $k>0$, then $I^\Pc_\gamma$ has infinite
index, and so does $I^{\Pc,\Qc}_\gamma$:
\begin{prop}\label{thm:indinfinito}
If $k>0$, then $I^\Pc_\gamma$ has infinite index in $\mathcal H_\gamma^\Pc$.
If $k\ge2$ or if $k=1$ and $\g(\dot\gamma,\dot\gamma)>0$, then $I^\Pc_\gamma$
has infinite index in the space of all vector fields in $\mathcal H_\gamma^\Pc$
that are everywhere orthogonal to $\dot\gamma$.
\end{prop}
\begin{proof}
If $Y$ is a Jacobi field along $\gamma$ and $f:[a,b]\to\R$ is a smooth function
vanishing at the endpoints, it is easily computed:
\begin{equation}\label{eq:formjaobi}
I_\gamma^\Pc(fY,fY)=\int_a^b\Big[{f'}^2\g(Y,Y)+\frac{\mathrm
d}{\mathrm dt}\,
\left(f^2\g(Y',Y)\right)\Big]\;\mathrm dt=\int_a^b{f'}^2\g(Y,Y)\;\mathrm dt.
\end{equation}
Let $t_0\in\left]a,b\right[$; if $k>0$, then we can find a Jacobi field $Y$ with
$\g(Y,Y)<0$ in a neighborhood $V$ of $t_0$. If $k\ge2$ or if $k=1$ and
$\g(\dot\gamma,\dot\gamma)>0$, then the field $Y$ can also be chosen orthogonal
to $\dot\gamma$ everywhere. From \eqref{eq:formjaobi}, it follows that $I_\gamma^\Pc$
is negative definite in the space of fields $fY$, where $f$ is supported in $V$.
\end{proof}
Obviously, the result of Proposition~\ref{thm:indinfinito} holds for
the bilinear form $I^{\Pc,\Qc}_\gamma$.
\end{section}
\begin{section}{The Maslov Index}\label{sec:maslov}
In this section we present some techniques of symplectic
spaces and we discuss the notion of {\em Maslov index\/}
that will be used to define an integer valued invariant
for semi-Riemannian geodesics.

\subsection{The Maslov index of a curve of Lagrangians}
\label{sub:Maslovcurve}
Let $(V,\omega)$ be a finite dimensional symplectic space, i.e.,
$V$ is a $2n$-dimensional real vector space and $\omega$ is a
nondegenerate skew-symmetric bilinear form in $V$. A subspace
$L$ of $V$ is {\em Lagrangian\/} if $\mathrm{dim}(L)=n$ and
$\omega$ vanishes on $L\times L$. The set $\Lambda$ of all Lagrangian
subspaces of $V$ is called the {\em Lagrangian Grassmannian\/} of $(V,\omega)$;
$\Lambda$ is a compact, connected
real analytic $\frac12n(n+1)$-dimensional embedded submanifold of the Grassmannian
$G_n(V)$ of all $n$-dimensional subspaces of $V$. We will use several well known facts
about the geometry  of the Lagrangian Grassmannian of a symplectic space (see for
instance
\cite{Ar, Duis, MPT});  in particular, we will make
full use of the notations and of the results proven in Reference~\cite{MPT}.

For our purposes, we need the following description of an atlas
of charts on $\Lambda$. Given a pair $L_0,L_1$ of complementary Lagrangian
subspaces of $V$, i.e., $V=L_0\oplus L_1$, we define an  isomorphism
$\Iota{L_0}{L_1}:L_1\to L_0^*$ by:
\begin{equation}\label{eq:defIota}
\Iota{L_0}{L_1}(v)=\omega(v,\cdot)\vert_{L_0},\quad
v\in L_1.
\end{equation}
We observe that, by the anti-symmetry of $\omega$, the following identity
holds:
\begin{equation}\label{eq:idIota}
\Iota{L_0}{L_1}=-(\Iota{L_1}{L_0})^*.
\end{equation}
Let $L\in\Lambda$
be fixed; we define the following subsets of $\Lambda$:
\begin{equation}\label{eq:defLambdak}
\Lambda_k(L)=\Big\{L'\in\Lambda:\mathrm{dim}(L'\cap L)=k\Big\},\quad k=0,\ldots,n.
\end{equation}
Each $\Lambda_k(L)$ is a connected embedded real analytic submanifold of $\Lambda$
having codimension $\frac12k(k+1)$ in $\Lambda$; $\Lambda_0(L)$ is a dense
open subset of $\Lambda$, while its complementary set:
\begin{equation}\label{eq:Lambda>=1}
\Lge(L)=\bigcup_{k=1}^n\Lambda_k(L)
\end{equation}
is {\em not\/} a regular submanifold of $\Lambda$, but only an {\em analytic subset}. Its
regular part is given by $\Lambda_1(L)$, which is a dense open subset of $\Lge(L)$.

Given  a pair $L_0,L_1$ of complementary Lagrangians in $V$, it is defined
a chart \[\phi_{L_0,L_1}:\Lambda_0(L_1)\to\Bsym(L_0,\R),\] 
where $\Bsym(L_0,\R)$ is the vector space of symmetric bilinear forms on
$L_0$. For $L\in\Lambda_0(L_1)$, we have:
\begin{equation}\label{eq:defphiLoL1}
\phi_{L_0,L_1}(L)=\Iota{L_0}{L_1}\circ T,
\end{equation}
where $T:L_0\to L_1$ is the unique linear map whose graph in $L_0\oplus L_1=V$
is equal to $L$. 
In equality \eqref{eq:defphiLoL1} we are identifying a linear map
$L_0\to L_0^*$ with a bilinear form on $L_0$; such identifications of
linear maps from a space to its dual and bilinear forms
on the space will be used throughout in the rest of the section.

Observe that, given $L\in\Lambda_0(L_1)$, the bilinear form
$\phi_{L_0,L_1}(L)$ is nondegenerate (i.e., the corresponding linear
map $L_0\to L_0^*$ is invertible) if and only if $L\in\Lambda_0(L_0)$.

The map $\phi_{L_0,L_1}$ defined in \eqref{eq:defphiLoL1} is a diffeomorphism, and it
follows in particular that $\Lambda_0(L_1)$ is contractible for all $L_1\in\Lambda$.
The Lagrangian Grassmannian $\Lambda$ is diffeomorphic to the homogeneous space
$U(n)/O(n)$ (\cite[Proposition~3.2.5]{MPT}), and using such diffeomorphism
one computes the fundamental group $\pi_1(\Lambda)\simeq\Z$
(\cite[Corollary~4.1.2]{MPT}).  It follows that the first singular homology
group $H_1(\Lambda;\Z)$ is also isomorphic to $\Z$; for a given Lagrangian
$L_0\in\Lambda$, since $\Lambda_0(L_0)$ is contractible, we compute the
first relative homology group of the pair $(\Lambda,\Lambda_0(L_0))$ as:
\begin{equation}\label{eq:relhomology}
H_1(\Lambda,\Lambda_0(L_0);\Z)\simeq \Z.
\end{equation}
The choice of the above isomorphism is related to the choice
of a {\em transverse orientation\/} of $\Lambda_1(L_0)$ in $\Lambda$, 
which is canonically associated to the symplectic form (\cite[Proposition~3.2.10]{MPT}). 
Every continuous curve $l$ in
$\Lambda$ with endpoints in $\Lambda_0(L_0)$ defines an element in
$H_1(\Lambda,\Lambda_0(L_0);\Z)$, and we denote by 
\begin{equation}\label{eq:maslovcurve}
\mu_{L_0}(l)\in\Z
\end{equation}
the integer number corresponding to the homology class of $l$ by
the isomorphism \eqref{eq:relhomology}. This number, which is additive
by concatenation and invariant by homotopies with endpoints in
$\Lambda_0(L_0)$,  can be interpreted as an {\em intersection number\/} 
of the curve $l$ with $\Lge(L_0)$.
\begin{defin}\label{thm:defmaslovcurve}
Given a continuous curve $l:[a,b]\to\Lambda$ with $l(a),l(b)\in\Lambda_0(L_0)$,
the integer number
$\mu_{L_0}(l)$ of \eqref{eq:maslovcurve} is called the
{\em Maslov index\/} of $l$ relative to $L_0$.
\end{defin}
The Maslov index of a continuous curve in $\Lambda$ can be computed
in terms of the coordinate charts $\phi_{L_0,L_1}$:
\begin{prop}\label{thm:calcMaslov}
Let $L_0\in\Lambda$ and let $l:[a,b]\to\Lambda$ be any continuous curve 
with endpoints in $\Lambda_0(L_0)$. Suppose that there exists a Lagrangian
subspace $L_1$ complementary to $L_0$ such that the image of $l$ is entirely
contained in the domain $\Lambda_0(L_1)$ of the chart $\phi_{L_0,L_1}$.
Then, the Maslov index $\mu_{L_0}(l)$ is given by:
\begin{equation}\label{eq:calcMaslov}
\mu_{L_0}(l)=n_+\Big(\phi_{L_0,L_1}\big(l(b)\big)\Big)-n_+\Big(\phi_{L_0,L_1}
\big(l(a)\big)\Big).
\end{equation}
\end{prop}
\begin{proof}
See \cite[Proposition~4.3.1]{MPT}
\end{proof}
To our purposes, we will need to extend the result of Proposition~\ref{thm:calcMaslov}
to the case that the image of the curve $l$ fails to be contained in
the domain of the chart $\phi_{L_0,L_1}$ at an isolated instant $t_0\in\left]a,b\right[$.
We need first  a technical Lemma:
\begin{lem}\label{thm:tech}
Let $L,L_*,L_0,L_1$ be four Lagrangian subspaces of $V$, with $L_0$ and $L_1$
complementary to each other, $L$ complementary to $L_0$ and with $L_*$ 
complementary to both $L$ and $L_0$. Then,
\begin{equation}\label{eq:tech}
\phi_{L_1,L_0}(L_*)-\phi_{L_1,L_0}(L)=(\Iota{L_0}{L_1})^*\circ\phi_{L_0,L_*}(L)^{-1}
\circ\Iota{L_0}{L_1}.
\end{equation}
\end{lem}
\begin{proof}
Let $T,S:L_1\to L_0$ be linear maps whose graphs in $V=L_1\oplus L_0$ are equal to 
$L_*$ and $L$ respectively; moreover let $U:L_0\to L_*$ be the linear
map whose graph in $V=L_0\oplus L_*$ is $L$. Observe that $U$ is invertible;
it is easily computed:
\begin{equation}\label{eq:1coisa}
Sv=U^{-1}(v+Tv)+Tv,\quad\forall\, v\in L_1.
\end{equation}
From \eqref{eq:defphiLoL1}, we have:
\begin{equation}\label{eq:2coisa}
\phi_{L_0,L_*}(L)=\Iota{L_0}{L_*}\circ U,\quad
\phi_{L_1,L_0}(L)=\Iota{L_1}{L_0}\circ S,\quad
\phi_{L_1,L_0}(L_*)=\Iota{L_1}{L_0}\circ T.
\end{equation}
From \eqref{eq:defIota}, we compute:
\begin{equation}\label{eq:3coisa}
(\Iota{L_0}{L_*})^{-1}\circ\Iota{L_0}{L_1}(v)=v+Tv,\quad \forall\,v\in L_1.
\end{equation}
Using \eqref{eq:1coisa}, \eqref{eq:2coisa} and \eqref{eq:3coisa}, it follows
\begin{equation}\label{eq:4coisa}
\phi_{L_1,L_0}(L)-\phi_{L_1,L_0}(L_*)=\Iota{L_1}{L_0}\circ\phi_{L_0,L_*}(L)^{-1}
\circ\Iota{L_0}{L_1}.
\end{equation}
The conclusion follows from \eqref{eq:idIota} and \eqref{eq:4coisa}.
\end{proof}
\begin{cor}\label{thm:corn+}
Under the assumptions of Lemma~\ref{thm:tech}, we have:
\begin{equation}\label{eq:corn+}
n_+\big(\phi_{L_1,L_0}(L_*)-\phi_{L_1,L_0}(L)\big)=
n_+\big(\phi_{L_0,L_*}(L)\big).
\end{equation}
In addition, $\phi_{L_1,L_0}(L_*)-\phi_{L_1,L_0}(L)$ is nondegenerate.
\end{cor}
\begin{proof}
It follows immediately from \eqref{eq:tech}, considering that:
\[n_+\big((\Iota{L_0}{L_1})^*\circ\phi_{L_0,L_*}(L)^{-1}
\circ\Iota{L_0}{L_1}\big)=n_+\big(\phi_{L_0,L_*}(L)^{-1}\big)=
n_+\big(\phi_{L_0,L_*}(L)\big).\qedhere\]
\end{proof}
We can now prove the aimed extension of Proposition~\ref{thm:calcMaslov}:
\begin{prop}\label{thm:calcMaslov2}
Let $L_0,L_1\in\Lambda$ be given, with $L_0\cap L_1=\{0\}$, and let $l:[a,b]\to\Lambda$
be a continuous curve  such that $l(t)\in\Lambda_0(L_0)$ except possibly for
$t=t_0\in\left]a,b\right[$. Let
$L_*\in\Lambda$ be complementary to both $l(t_0)$ and $L_0$; then,  for $\varepsilon>0$
sufficiently small, we have:
\begin{equation}\label{eq:calcMaslov2}
\mu_{L_0}(l)=n_-\Big(\phi_{L_1,L_0}\big(l(t_0+\varepsilon)\big)-\phi_{L_1,L_0}\big(L_*
\big)\Big)-n_-\Big(
\phi_{L_1,L_0}\big(l(t_0-\varepsilon)\big)-\phi_{L_1,L_0}\big(L_*\big)\Big).
\end{equation}
\end{prop}
\begin{proof}
Let $\varepsilon>0$ be small enough so that $l(t)\in\Lambda_0(L_*)$
for all $t\in[t_0-\varepsilon,t_0+\varepsilon]$; since $t_0$ is the unique
instant where $l$ passes through $\Lge(L_0)$, then $\mu_{L_0}(l)=\mu_{L_0}
(l\vert_{[t_0-\varepsilon,t_0+\varepsilon]})$. Using Proposition~\ref{thm:calcMaslov},
we obtain 
\[\mu_{L_0}(l)=n_+\Big(\phi_{L_0,L_*}\big(l(t_0+\varepsilon)\big)\Big)-
n_+\Big(\phi_{L_0,L_*}\big(l(t_0-\varepsilon)\big)\Big).\]
The conclusion follows by applying twice Corollary~\ref{thm:corn+} to the above
equation, once by taking $L=l(t_0+\varepsilon)$ and again by taking
$L=l(t_0-\varepsilon)$.
\end{proof}
\subsection{The Maslov index of a semi-Riemannian geodesic}
\label{sub:Maslovgeo}
We now consider a semi-Rie\-mann\-ian setup as in Section~\ref{sec:semiriemannian},
consisting of a semi-Riemannian manifold $(\M,\g)$, a nondegenerate
smooth submanifold $\Pc$ of $\M$ and a geodesic $\gamma:[a,b]\to\M$ starting
orthogonally to $\Pc$. We start by observing that, for $J_1,J_2\in\J$, we have:
\begin{equation}\label{eq:simetria}
\g\big(J_1'(t),J_2(t)\big)=\g\big(J_1(t),J_2'(t)\big),\quad\forall\,t\in[a,b].
\end{equation}
We choose a parallel trivialization of the tangent bundle $T\M$ along
$\gamma$; we may then identify vector fields along $\gamma$ with curves
in $\R^n$ and the metric tensor $\g$ along $\gamma$ with a fixed
nondegenerate symmetric bilinear form $g$ in $\R^n$. The space
$\J$ will then correspond to a space $\Jb$ of smooth curves in $\R^n$.

Let us consider the canonical symplectic structure $\omega$ on the vector space
$V=\R^n\oplus{\R^n}^*$ given by:
\begin{equation}\label{eq:defomega}
\omega\big((v_1,\alpha_1),(v_2,\alpha_2)\big)=\alpha_2(v_1)-\alpha_1(v_2).
\end{equation}
For all $t\in[a,b]$, we define an $n$-dimensional subspace $\ell(t)\subset V$ by:
\begin{equation}\label{eq:defellt}
\ell(t)=\Big\{\big(J(t),gJ'(t)\big):J\in\Jb\Big\};
\end{equation}
here $g$ is thought of as a linear map from $\R^n$ to ${\R^n}^*$.
By \eqref{eq:simetria}, $\ell(t)$ is a Lagrangian subspace of $(V,\omega)$
for all $t\in[a,b]$, and we therefore obtain a smooth curve $\ell:[a,b]\to\Lambda$.

We fix the following  Lagrangian subspace $L_0$ of $V$:
\begin{equation}\label{eq:L0}
L_0=\{0\}\oplus{\R^n}^*.
\end{equation}
Observe that, for $t\in\left[a,b\right]$, $\ell(t)\in\Lge(L_0)$ if and only if
$\gamma(t)$ is a $\Pc$-focal point; moreover, the multiplicity of
$\gamma(t)$ coincides with the dimension of $\ell(t)\cap L_0$.
In particular, if $\gamma(b)$ is not a $\Pc$-focal point, then the curve
$\ell$ has final endpoint in $\Lambda_0(L_0)$. On the other hand, 
$\ell(a)\in\Lge(L_0)$; however, since there are no $\Pc$-focal points near
$a$ (Proposition~\ref{thm:isolPfocalpts}), in order to define the
Maslov index of the geodesic $\gamma$ we can consider a restriction
$\ell\vert_{[a+\varepsilon,b]}$ with $\varepsilon>0$ small.
\begin{defin}\label{thm:defMaslovgeo}
Suppose that $\gamma(b)$ is not $\Pc$-focal. The {\em Maslov index\/}
$\imaslov(\gamma)$ of the geodesic $\gamma$ is defined as:
\begin{equation}\label{eq:defMaslovgeo}
\imaslov(\gamma)=\mu_{L_0}\big(\ell\vert_{[a+\varepsilon,b]}\big),
\end{equation}
where $\varepsilon>0$ is chosen such that $\gamma(t)$ is not
$\Pc$-focal for $t\in\left]a,a+\varepsilon\right]$.
\end{defin}
Clearly, the right hand side of \eqref{eq:defMaslovgeo} does not
depend on the choice of $\varepsilon$; moreover, in order to make
rigorous the above definition we need the following:
\begin{prop}\label{thm:naodep}
The term on the right hand side of equality \eqref{eq:defMaslovgeo}
does not depend on the choice of a parallel trivialization
of $T\M$ along $\gamma$.
\end{prop} 
\begin{proof}
If $\tilde\ell:[a,b]\to\Lambda$ is the curve of Lagrangians
corresponding to a different choice of a parallel trivialization of $T\M$
along $\gamma$,
then the relation between $\ell$ and $\tilde\ell$ is given by:
\[\tilde\ell=\sigma\circ\ell,\]
where $\sigma:V\to V$ is a fixed symplectomorphism that preserves
$L_0$. Namely, $\sigma$ is given by:
\[\sigma(v,\alpha)=\big(s(v),{s^*}^{-1}(\alpha)\big),\]
where $s:\R^n\to\R^n$ is the isomorphism that relates the two trivializations.
The conclusion follows from the fact that composition with a fixed symplectomorphism
that preserves $L_0$ induces the identity in the relative homology group
$H_1(\Lambda,\Lambda_0(L_0))$ (\cite[Remark~4.2.1]{MPT}).
\end{proof}
We have the following relation between the Maslov index and the focal
index of a semi-Riemannian geodesic:
\begin{prop}\label{thm:imaslovifoc}
Suppose that $\gamma(b)$ is not $\Pc$-focal and that all the $\Pc$-focal
points are nondegenerate. Then,
\begin{equation}\label{eq:imaslovifoc}
\imaslov(\gamma)=\ifoc(\gamma).
\end{equation}
\end{prop}
\begin{proof}
See \cite[Theorem~5.1.2]{MPT}.
\end{proof}
We remark that the thesis of Proposition~\ref{thm:imaslovifoc} is false
without the nondegeneracy assumption on the focal points, even for real
analytic manifolds (see
\cite[Subsection~7.4]{MPT}). 

It is easy to see that, due to its topological nature, the Maslov
index is invariant by uniformly small perturbations of the data
of the geometric problem; on the other hand, the focal index
is unstable. The stability property is a first indication that
the Maslov index is the correct generalization of the notion
of geometric index to semi-Riemannian geometry. 
\end{section}
\begin{section}{Abstract Results of Functional Analysis}\label{sec:abstract}
The goal of this section is to provide a method of computing the
change of index of a smooth family of symmetric bilinear forms
on a Hilbert space; we will use the notations and several results 
from \cite{GMPT} that will be restated for the reader's convenience.
All Hilbert spaces of the Section will be assumed {\em real}.

Given Hilbert spaces $\mathcal H,\mathcal H'$, we will denote
by $\mathcal L(\mathcal H,\mathcal H')$ the space of
bounded linear operators from $\mathcal H$ to $\mathcal H'$; by $\mathcal L(\mathcal H)$
we will mean $\mathcal L(\mathcal H,\mathcal H)$. By $\Bsym(\mathcal H,\R)$ we will now
mean the set of symmetric  {\em bounded\/} bilinear forms on
$\mathcal H$. 
Let $\langle\cdot,\cdot\rangle$ be a Hilbert space inner product on $\mathcal H$;
to any bounded bilinear form $B:\mathcal H\times \mathcal H\to\R$ by Riesz's theorem
there corresponds a bounded linear operator $T_B:\mathcal H\to\mathcal H$, which
is related to $B$ by:
\begin{equation}\label{eq:defBT}
B(x,y)=\langle T_B(x),y\rangle,\quad\forall\,x,y\in\mathcal H.
\end{equation}
We say that $T_B$ is the linear operator that {\em represents\/}  $B$ with respect
to the inner product $\langle\cdot,\cdot\rangle$. Clearly, $B$ is symmetric
if and only if $T_B$ is self-adjoint. We say that $B$ is {\em nondegenerate\/}
if $T_B$ is injective; $B$ will be said to be {\em strongly nondegenerate\/}
if $T_B$ is an isomorphism. If $T_B$ is a {\em Fredholm operator\/} of index
$0$ (for instance if $T_B$ is a compact perturbation of an isomorphism), then
 $B$ is nondegenerate if and only if it is strongly
nondegenerate. Observe that  strong nondegeneracy is stable by small
perturbations, since the set of isomorphisms of $\mathcal H$ is open
in $\mathcal L(\mathcal H)$.\smallskip

We will consider $1$-parameter families of bilinear forms defined
on a variable domain, and we need the following notion of
$C^1$-family of closed subspace of a Hilbert space:
\begin{defin}\label{thm:defC1subspaces}
Let $\mathcal H$ be a Hilbert space, $I\subset\R$ an interval and $\{\mathcal D_t\}_{t\in I}$
be a family of closed subspaces of $\mathcal H$. We say that $\{\mathcal D_t\}_{t\in I}$
is a {\em $C^1$-family\/} of subspaces if for all $t_0\in I$ there exists a $C^1$-curve
$\alpha:\left]t_0-\varepsilon,t_0+\varepsilon\right[\,\cap\, I\to \mathcal L(\mathcal
H)$  and a closed subspace $\overline{\mathcal D}\subset\mathcal H$
such that $\alpha(t)$ is an isomorphism  and $\alpha(t)(\mathcal D_t)=\overline{\mathcal D}$
for all~$t$.
\end{defin}
We have the following criterion to establish the regularity of
a family of closed subspaces:
\begin{lem}\label{thm:produce}
Let $I\subset\R$ be an interval, $\mathcal H,\tilde{\mathcal H}$ be Hilbert
spaces and $F:I\to\mathcal L(\mathcal H,\tilde{\mathcal H})$ be a $C^1$-map
such that each $F(t)$ is surjective. Then, the family $\mathcal D_t=\mathrm{Ker}(F(t))$
is a $C^1$-family of closed subspaces of $\mathcal H$.
\end{lem}
\begin{proof}
See~\cite[Lemma~2.9]{GMPT}.
\end{proof}
The next Proposition, also proven in \cite{GMPT}, gives a method
for computing the change of the index of a smooth family of bilinear forms
that are represented by a {\em compact perturbation of a positive
isomorphism}. Recall that a self-adjoint linear operator $T$ in $\mathcal H$
is a compact perturbation of a positive (negative) isomorphism of $\mathcal H$ if it is
of the form $L+K$, where $L$ is a self-adjoint positive (negative) isomorphism
of $\mathcal H$
 and $K$ is a compact self-adjoint operator on $\mathcal H$.
\begin{prop}\label{thm:HSelementary}
Let $\mathcal H$ be a  Hilbert space with inner product
$\langle\cdot,\cdot\rangle$, and let $B:[t_0,t_0+r]\to\Bsym(\mathcal H,\R)$, $r>0$,
be a map of class $C^1$.  Let $\{\mathcal D_t\}_{t\in[t_0,t_0+r]}$ be a $C^1$-family of
closed subspaces of $\mathcal H$, and denote by $\overline B(t)$ the restriction of
$B(t)$ to $\mathcal D_t\times\mathcal D_t$.
Assume that the following three hypotheses are satisfied:
\begin{enumerate}
\item\label{itm:hp2.5.1} $\overline B(t_0)$ is represented by a
compact perturbation of a positive isomorphism of $\mathcal D_{t_0}$;
\item\label{itm:hp2.5.2} the restriction $\widetilde B$ of the derivative $B'(t_0)$ to 
$\mathrm{Ker}(\overline B(t_0))\times \mathrm{Ker}(\overline B(t_0))$ 
is nondegenerate;
\item\label{itm:hp2.5.3} $\mathrm{Ker}\big(\overline B(t_0)\big)\subseteq
\mathrm{Ker}(B(t_0))$.
\end{enumerate}
Then, for
$t>t_0$ sufficiently close to $t_0$, $\overline B(t)$ is nondegenerate, and we have:
\begin{equation}\label{eq:changen-}
n_-\big(\overline B(t)\big)=n_-\big(\overline B(t_0)\big)+n_-(\widetilde B),
\end{equation}
all the terms of the above equality being finite natural numbers.
\end{prop}
\begin{proof}
See \cite[Proposition~2.5]{GMPT}.
\end{proof}
\begin{rem}\label{thm:remcasoparticular}
Observe that, by Proposition~\ref{thm:HSelementary}, if $\overline B(t_0)$ is
nondegenerate on $\mathcal D_{t_0}$, then $n_-(\overline B(t))$ is constant
for $t$ near $t_0$. Actually, we have the following stronger continuity property
for the positive and the negative type numbers of symmetric bilinear forms.
If $B_n\to B$ in $\Bsym(\mathcal H,\R)$,
$\mathcal D_n$ converges\footnote{in the sense that $F_n\to F$ in $\mathcal L(\mathcal
H,\tilde{\mathcal H})$, where $\tilde{\mathcal H}$ is any Hilbert space,
$F$ is surjective and $\mathcal D_n=\mathrm{Ker}(F_n)$, $\mathcal D=\mathrm{Ker}(F)$.} to
$\mathcal D$, if
$B\vert_{\mathcal D}$ is nondegenerate and it is represented by a compact
perturbation of a positive (resp., negative) isomorphism of $\mathcal D$, then
for $n$ sufficiently large, it is $n_-(B_n\vert_{\mathcal D_n})=n_-(B\vert_{\mathcal D})$
(resp., $n_+(B_n\vert_{\mathcal D_n})=n_+(B\vert_{\mathcal D})$).
\end{rem}
For the purposes of this article, we need an extension of the result
of Proposition~\ref{thm:HSelementary} that holds in the more general
situation in which  hypothesis~\eqref{itm:hp2.5.3} is not satisfied.
To this aim, we need to define a notion of {\em derivative\/} of
the family $B(t)$ that takes into consideration the variation of the
domain $\mathcal D_t$. 

\begin{defin}\label{thm:defderivada}
Let $\mathcal H$ be a  Hilbert space and let $B:[t_0,t_0+r]\to \Bsym(\mathcal
H,\R)$, $r>0$, be a map of class $C^1$.  Let $\{\mathcal D_t\}_{t\in[t_0,t_0+r]}$ be a
$C^1$-family of closed subspaces of $\mathcal H$, and denote by $\overline B(t)$ the
restriction of $B(t)$ to $\mathcal D_t\times\mathcal D_t$. We define the symmetric
bilinear form
$\overline B'(t_0)$ in $\mathrm{Ker}(\overline B(t_0))$ by:
\begin{equation}\label{eq:defderivada}
\begin{split}
\overline B'&(t_0)(v,w)=\frac{\mathrm d}{\mathrm
dt}\,B(t)\big(v(t),w(t)\big)\Big\vert_{t=t_0}=\\
&=B'(t_0)(v,w)+B(t_0)\big(v'(t_0),w\big)+B(t_0)\big(v,w'(t_0)\big),\quad\forall\,v,w\in
\mathrm{Ker}(\overline B(t_0)),
\end{split}
\end{equation}
where $v(t)$ and $w(t)$ are $C^1$-curves in $\mathcal H$ with 
$v(t_0)=v$, $w(t_0)=w$, $v(t)\in\mathcal D_t$ and $w(t)\in\mathcal D_t$ for
all $t$.
\end{defin}
\begin{rem}\label{thm:remindepchoice}
Note that  formula \eqref{eq:defderivada} defines $\overline B'(t_0)(v,w)$
independently of the extensions $v(t)$ and $w(t)$ chosen. To see this,
simply observe that the classes of the derivatives $v'(t_0)$ and $w'(t_0)$ 
modulo $\mathcal D_{t_0}$ are independent of the extensions.
\end{rem}

This is the aimed extension of Proposition~\ref{thm:HSelementary}:
\begin{prop}\label{thm:notsoelementary}
Let $\mathcal H$ be a  Hilbert space with inner product
$\langle\cdot,\cdot\rangle$, and let $B:[t_0,t_0+r]\to\Bsym(\mathcal H,\R)$, $r>0$,
be a map of class $C^1$.  Let $\{\mathcal D_t\}_{t\in[t_0,t_0+r]}$ be a $C^1$-family of
closed subspaces of $\mathcal H$, and denote by $\overline B(t)$ the restriction of
$B(t)$ to $\mathcal D_t\times\mathcal D_t$.
Assume that the following two hypotheses are satisfied:
\begin{enumerate}
\item\label{itm:hp2.5.1b} $\overline B(t_0)$ is represented by a
compact perturbation of a positive isomorphism of $\mathcal D_{t_0}$;
\item\label{itm:hp2.5.2b} the symmetric bilinear form $\overline B'(t_0)$
is nondegenerate.
\end{enumerate}
Then, for
$t>t_0$ sufficiently close to $t_0$, $\overline B(t)$ is nondegenerate, and we have:
\begin{equation}\label{eq:changen-b} 
n_-\big(\overline B(t)\big)=n_-\big(\overline B(t_0)\big)+n_-\big(\overline B'(t_0)\big),
\end{equation}
all the terms of the above equality being finite natural numbers.
\end{prop}
\begin{proof}
By possibly passing to a smaller $r$, we can assume the existence of
a $C^1$-curve $\alpha(t)$ of isomorphisms of $\mathcal H$ such that $\alpha(t)$
carries $\mathcal D_t$ to a fixed subspace $\overline{\mathcal D}$ of $\mathcal H$.
Define  $C(t)=B(t)\big(\alpha(t)^{-1}\cdot,\alpha(t)^{-1}\cdot\big)$ as a bilinear form
on the fixed space $\overline{\mathcal D}$. 
Then, $C(t)$ is a push-forward of $\overline B(t)$ and the restriction
of $C'(t_0)$ to $\mathrm{Ker}\big(C(t_0)\big)$ 
is a push-forward of $\overline B'(t_0)$. The conclusion follows by applying
Proposition~\ref{thm:HSelementary}\/ to the curve $C(t)$ in $\Bsym(\overline{\mathcal
D};\R)$.
\end{proof}

\begin{cor}\label{thm:corelementary}
Let $B:[t_0-r,t_0+r]\to\Bsym(\mathcal H,\R)$ and $\{\mathcal D_t\}_{t\in[t_0-r,t_0+r]}$
satisfy the same hypotheses of Proposition~\ref{thm:notsoelementary}. Then, in the
notations of Proposition~\ref{thm:notsoelementary},
for $\varepsilon>0$ small enough, we have:
\begin{equation}\label{eq:doublechangen-}
n_-\big(\overline B(t_0-\varepsilon)\big)-n_-\big(\overline B(t_0+\varepsilon)\big)=\sgn
\big(\overline B'(t_0)\big).
\end{equation}
\end{cor}
\begin{proof}
Use Proposition~\ref{thm:notsoelementary} twice, once to $B\vert_{[t_0,t_0+r]}$
and once to a backwards  reparameterization of $B\vert_{[t_0-r,t_0]}$.
\end{proof}

\end{section}
\begin{section}{The Morse Index Theorem}\label{sec:index}
In this section we go back to the geometrical setup
of Section~\ref{sec:semiriemannian} and we state and prove an extension of the
Morse Index Theorem for geodesics in semi-Riemannian manifolds 
with metric tensor of arbitrary index. As we have seen in
Proposition~\ref{thm:indinfinito}, if $(\M,\g)$ is not Riemannian then the index of
$I^\Pc_\gamma$ is always infinite. However, we show that it is possible to split
the Hilbert space of all variations of a given geodesic into two
subspaces such that  $I_\gamma^\Pc$ has finite index on the first and
finite coindex on the second. 

The definition of these spaces of variations depend on the choice of 
a distribution of maximal negative subspaces along the geodesic $\gamma$:
\begin{defin}\label{thm:defmaxnegdistr}
We say that a  family of subspaces $\mathcal D_t\subset 
T_{\gamma(t)}\M$, $t\in[a,b]$, along the geodesic $\gamma$
is {\em smooth\/} if there exist a family $Y_1,\ldots,Y_r$ of smooth
vector fields along $\gamma$ which forms a pointwise basis for $\mathcal D$;
such a family $Y_1,\ldots,Y_r$ is called a {\em frame\/} for $\mathcal D$.
A {\em maximal negative distribution\/} along  $\gamma$
is a smooth family of $k$-dimensional subspaces $\mathcal D$ along $\gamma$
such that $\g$ is negative definite on $\mathcal D_t$ for all~$t$
(recall \eqref{eq:indexg}). 
\end{defin}

Obviously, maximal negative distributions along any geodesic always exist;
for instance, one can obtain such distributions by considering
the parallel transport of any maximal subspace of $T_{\gamma(a)}\M$
on which $\g$ is negative definite.

Given a maximal negative distribution $\mathcal D$ along $\gamma$,
we define the following closed subspaces of $\mathcal H_\gamma^\Pc$:
\begin{equation}\label{eq:defKS}
\begin{split}
\KgD=\Big\{v\in\mathcal H_\gamma^\Pc:&\;\g(v',Y_i)\ \text{is of Sobolev regularity
$H^1$, and}\\  &\;\g(v',Y_i)'=\g(v',Y_i')+\g\big(\mathcal
R(\dot\gamma,v)\,\dot\gamma,Y_i\big),\ i=1,\ldots,k\Big\}
\\
\SgD=\Big\{v\in\mathcal H_\gamma^\Pc:&\;v(a)=0,\ v(t)\in\mathcal D_t,\ \
\forall\,t\in[a,b]\Big\},
\end{split}
\end{equation} 
where $Y_1,\ldots,Y_k$ is a frame for $\mathcal D$. It is easy to check that
the space $\KgD$ does {\em not\/} actually depend
on the choice of the frame $Y_1,\ldots,Y_k$. 
The space $\KgD$ can be roughly described as the space of vector fields
along $\gamma$ that are ``Jacobi in the directions of $\mathcal D$'';
observe indeed that if $v\in\mathcal H^\Pc_\gamma$
is a vector field of class $C^2$, then $v\in\KgD$ if and only if:
\[v''-\mathcal R(\dot\gamma,v)\,\dot\gamma\in\mathcal D^\perp.\]

We are ready to state the main result of the paper:

\begin{teo}[Semi-Riemannian Morse Index Theorem]\label{thm:indexth}
Let $(\M,\g)$ be a semi-Riemann\-ian manifold, $\Pc$ a smooth submanifold
of $\M$, $\gamma:[a,b]\to\M$ a geodesic such that:
\begin{itemize}
\item $\gamma(a)\in\Pc$ and $\dot\gamma(a)\in T_{\gamma(a)}\Pc^\perp$;
\item $\Pc$ is nondegenerate at $\gamma(a)$;
\item $\gamma(b)$ is not a $\Pc$-focal point.
\end{itemize}
Let $\mathcal D$ be a maximal negative distribution along $\gamma$; let $\KgD$
and $\SgD$ be the corresponding subspaces of $\mathcal H_\gamma^\Pc$ defined
in \eqref{eq:defKS}. Then,
\begin{equation}\label{eq:indexth}
\imaslov(\gamma)=n_-\left(I_\gamma^\Pc\big\vert_{\KgD}\right)-n_+
\left(I_\gamma^\Pc\big\vert_{\SgD}\right)
-n_-\left(\g\big\vert_{T_{\gamma(a)}\Pc}\right),
\end{equation}
where all the terms in the above formula are finite integer numbers.
\end{teo}

The proof of Theorem~\ref{thm:indexth} requires some work
and it is spread along the remaining subsections of this section.

The spaces $\KgD$ and $\SgD$ are $I^\Pc_\gamma$-orthogonal 
(Lemma~\ref{thm:KSorth}); moreover, under {\em generic\/} circumstances,
they are complementary in $\mathcal H_\gamma^\Pc$ (Corollary~\ref{thm:multfocinstred}
and Corollary~\ref{thm:qFsurjective}). An explicit formula to compute
the term $n_+\left(I_\gamma^\Pc\big\vert_{\SgD}\right)$ that appears
in \eqref{eq:indexth} is given in Corollary~\ref{thm:indIS}.

The last term of equality \eqref{eq:indexth} is the contribution
given by the initial submanifold $\Pc$; in the case of Riemannian
or causal Lorentzian geodesics, $\Pc$ is spacelike at $\gamma(a)$, and
therefore the last term of \eqref{eq:indexth} vanishes.
\smallskip

Let's take a closer look at some special examples to get
a better feeling of the result of Theorem~\ref{thm:indexth}.
\begin{example}\label{exa:riemannian}
If $(\M,\g)$ is Riemannian, then $\mathcal D=0$, the space $\KgD$ coincides
with $\mathcal H^\Pc_\gamma$, $\SgD=\{0\}$ and the Maslov index
of $\gamma$ is equal to the sum of the multiplicities of the $\Pc$-focal
points along $\gamma$. 
\end{example}
For simplicity, in our next example we will assume that the initial submanifolds
to the given geodesics reduce to a point; in this case we will omit
the subscripts and the superscripts $\mathcal P$ in our notation.
\begin{example}\label{exa:product}
Let $(\M_1,\g_1)$, $(\M_2,\g_2)$ be Riemannian manifolds and consider the product
$\M=\M_1\times\M_2$ endowed with the semi-Riemannian metric $\g=\g_1\oplus(-\g_2)$. It
is easily seen that $\gamma=(\gamma_1,\gamma_2):[a,b]\to\M$ is a geodesic iff $\gamma_1$
and $\gamma_2$ are geodesics; the space $\mathcal H_\gamma$ of vector fields along
$\gamma$ vanishing at both endpoints is identified with the direct sum $\mathcal
H_{\gamma_1}\oplus\mathcal H_{\gamma_2}$ where $\mathcal H_{\gamma_i}$ is the space of
vector fields along $\gamma_i$ vanishing at both endpoints, $i=1,2$. It is easily seen
that the index form $I_\gamma$ is given by $I_\gamma=I_{\gamma_1}\oplus(-I_{\gamma_2})$
where $I_{\gamma_i}$ is the index form corresponding to $\gamma_i$ in the {\em
Riemannian\/} manifold $(\M_i,\g_i)$, $i=1,2$; assuming that $\gamma(b)$ is not
conjugate to $\gamma(a)$ one easily sees that the spaces $\mathcal K_\gamma^{\mathcal
D}$ and $\mathcal S_\gamma^{\mathcal D}$ corresponding to the distribution $\mathcal
D_t=T_{\gamma_2(t)}\M_2\subset T_{\gamma(t)}\M$ are given respectively by $\mathcal
H_{\gamma_1}$ and
$\mathcal H_{\gamma_2}$. In this case, Theorem~\ref{thm:indexth} is an easy consequence
of the Riemannian Morse index theorem applied to each geodesic $\gamma_i$.
\end{example}
\begin{example}\label{exa:causal}
If $(\M,\g)$ is Lorentzian and $\gamma$ is timelike,
i.e., $\g(\dot\gamma,\dot\gamma)<0$, then we can consider the distribution
$\mathcal D$ spanned by $\dot\gamma$. In this case, $\KgD$ correspond
to the space of variational vector fields that are everywhere
orthogonal to $\dot\gamma$ and $n_+\left(I_\gamma^\Pc\big\vert_{\SgD}\right)=0$.
Also in this case, the Maslov index of $\gamma$ equals the sum of the multiplicities
of the $\Pc$-focal points along $\gamma$.
\end{example}

\begin{example}\label{exa:semi}
Suppose that, in the general semi-Riemannian case, 
we can find $Y_1,\ldots,Y_k$  Jacobi fields along $\gamma$, with
$k=n_-(\g)$, that form a frame for a $k$-dimensional distribution $\mathcal D$ on which
$\g$ is negative definite. If $\left(\g(Y_i',Y_j)\right)_{ij}$ is symmetric,
then also in this situation $n_+\left(I_\gamma^\Pc\big\vert_{\SgD}\right)=0$
(this will follow from Corollary~\ref{thm:corcriterio} ahead);
observe that we can always find such a family of Jacobi fields on
sufficiently small segments of a geodesic. 
In this context, the space $\KgD$ is given by the set of vectors fields $v\in \mathcal
H^\Pc_\gamma$ such that the quantities $\g(v',Y_i)-\g(v,Y_i')$ are constant for every
$i$. In the Lorentzian case, $k=1$ and this observation applies when
the geodesic $\gamma$ admits a timelike Jacobi field along it.
\end{example}

\begin{example}\label{exa:lie}
Suppose that $G$ is a $k$-dimensional Lie group acting on $\M$ by isometries
with no fixed points, or more in general, having only discrete
isotropy groups. Suppose that $\g$ is negative definite on the orbits of $G$.
If $\dot\gamma(a)$ is orthogonal to the orbit of the commutator subgroup
$[G,G]$ (for instance if $G$ is abelian), then we can consider the distribution 
$\mathcal D$ tangent to the orbits of $G$. Observe that $\mathcal D$ is generated
by $k$ linearly independent Killing vector fields $Y_1,\ldots,Y_k$
on $\M$, which therefore restrict to Jacobi fields along any geodesic.
Then, one falls into the case of Example~\ref{exa:semi}
by observing  that the symmetry of $\g(Y_i',Y_j)$ follows from
the orthogonality of $\dot\gamma(a)$ with the orbits of $[G,G]$:
\[\g(Y_i',Y_j)-\g(Y_i,Y_j')=-\g(\nabla_{Y_j}Y_i,\dot\gamma)+\g(\nabla_{Y_i}Y_j,\dot\gamma)=
\g([Y_i,Y_j],\dot\gamma)=0.\]
In this situation, the space $\KgD$ can be described
as the space of variational vector fields along $\gamma$ corresponding to
variations of $\gamma$ by curves that are {\em geodesics along $\mathcal D$}, i.e.,
whose second derivatives are orthogonal to $\mathcal D$.
\end{example}

\begin{example}\label{exa:n+0}
Another situation in which the term $n_+\left(I_\gamma^\Pc\big\vert_{\SgD}\right)$
vanishes occurs when the bilinear form $\g(\Rc(\dot\gamma,\cdot)\,\dot\gamma,\cdot)$ 
is negative semi-definite along the geodesic $\gamma$ and $\mathcal D$
is parallel (again, this will follow from Corollary~\ref{thm:corcriterio}).
\end{example}

We conclude this subsection by showing that Theorem~\ref{thm:indexth}
can be easily generalized to the case of geodesics with both endpoints variable.

In order to give a statement of this extension we 
need to introduce the following objects.

Assume that we are given a smooth submanifold $\Qc$ of $\M$ such that $\gamma(b)\in\Qc$
and $\dot\gamma(b)\in T_{\gamma(b)}\Qc^\perp$. In analogy with \eqref{eq:defKS},
we define the space $\KgDPQ$ by:
\begin{equation}\label{eq:defKgDPQ}
\begin{split}
\KgDPQ=\Big\{v\in\mathcal H_\gamma^{\Pc,\Qc}:&\;\g(v',Y_i)\ \text{is of Sobolev
regularity
$H^1$, and}\\  &\;\g(v',Y_i)'=\g(v',Y_i')+\g\big(\mathcal
R(\dot\gamma,v)\,\dot\gamma,Y_i\big),\ i=1,\ldots,k\Big\}
\end{split}
\end{equation}

Suppose that $\gamma(b)$ is not $\Pc$-focal;
let $S_\gamma$ be the linear endomorphism of
$T_{\gamma(b)}\M$ defined by:
\[S_\gamma\big(J(b)\big)=-J'(b),\]
for all $J\in\J$. Observe that the assumption of non focality
for $\gamma(b)$ implies that $\J\ni J\mapsto J(b)\in T_{\gamma(b)}\M$
is an isomorphism, and therefore $S_\gamma$ is well defined.
Observe also that, by \eqref{eq:simetria}, $S_\gamma$ is $g$-symmetric;
we denote by $S_\gamma$ also the corresponding symmetric bilinear form
on $T_{\gamma(b)}\M$, which is given by:
\[\phantom{\quad\forall\,J_1,J_2\in\J.}S_\gamma\big(J_1(b),J_2(b)\big)=-g\big(J_1(b),
J_2'(b)\big),\quad\forall\,J_1,J_2\in\J.\]

\begin{teo}\label{thm:indexPQ}
Under the hypotheses of Theorem~\ref{thm:indexth}, assume also that
we are given a smooth submanifold $\Qc$ of $\M$ such that $\gamma(b)\in\Qc$
and $\dot\gamma(b)\in T_{\gamma(b)}\Qc^\perp$. Then,
\begin{equation}\label{eq:indexPQ}
\begin{split}
\imaslov(\gamma)=&\;n_-\left(I_\gamma^{\Pc,\Qc}\big\vert_{\KgDPQ}\right)-n_+
\left(I_\gamma^{\Pc,\Qc}\big\vert_{\SgD}\right)
-n_-\left(\g\big\vert_{T_{\gamma(a)}\Pc}\right)\\&-n_-\left( \mathcal
S^\Qc_{\dot\gamma(b)}- S_\gamma \,\big\vert_{T_{\gamma(b)}\Qc}\right).
\end{split}
\end{equation}
\end{teo}
\begin{proof}
Recalling \eqref{eq:J}, we define:
\[\mathfrak J_\Qc=\Big\{J\in\mathfrak J:J(b)\in T_{\gamma(b)}\Qc\Big\}.\]
Since $\gamma(b)$ is not $\Pc$-focal and $\mathfrak J_\Qc\subset \Kt^{\mathcal
D}_{\gamma,\Pc,\Qc}$, it follows easily that 
\begin{equation}\label{eq:dirsum}
\KgDPQ=\KgD\oplus \mathfrak J_\Qc.
\end{equation}
Integration by parts in \eqref{eq:defIPQ} shows that the direct sum
in \eqref{eq:dirsum} is $I^{\Pc,\Qc}_\gamma$-orthogonal, hence
\begin{equation}\label{eq:soman-}
n_-\left(I^{\Pc,\Qc}_\gamma\big\vert_{\KgDPQ}\right)=n_-
\left(I^{\Pc,\Qc}_\gamma\big\vert_{\KgD}\right)+
n_-\left(I^{\Pc,\Qc}_\gamma\big\vert_{\mathfrak J_\Qc}\right).
\end{equation}
The restriction of $I^{\Pc,\Qc}_\gamma$ to $\KgD$ is obviously equal to the
restriction of $I_\gamma^\Pc$ to the same space, hence the first term on the right hand
side of equality
\eqref{eq:soman-} is computed in Theorem~\ref{thm:indexth}.

The conclusion follows by observing that the isomorphism
$\mathfrak J_\Qc\ni J\mapsto J(b)\in T_{\gamma(b)}\Qc$ carries
the restriction of $I_\gamma^{\Pc,\Qc}$ to $\mathcal
S^\Qc_{\dot\gamma(b)}- S_\gamma \,\big\vert_{T_{\gamma(b)}\Qc}$.
\end{proof}
Observe that the last term in equality \eqref{eq:indexPQ} is the contribution
of the final manifold $\Qc$; it already appears in the Riemannian
Morse Index Theorem for variable endpoints (\cite{Kal}).

We now pass to the proof of Theorem~\ref{thm:indexth}.
\subsection{Reduction to a Morse--Sturm system in $\mathbf{I\!R^n}$}\label{sub:reduction}
A Morse--Sturm system in $\R^n$ is a second order linear
differential system of the form:
\begin{equation}\label{eq:MS}
v''(t)=R(t)\,v(t),\quad t\in[a,b],\ v(t)\in\R^n,
\end{equation}
where $R(t)$ is a continuous map of linear endomorphisms of $\R^n$ that are
symmetric with respect to a fixed nondegenerate symmetric bilinear form
$g$ on $\R^n$.

Morse--Sturm systems arise from the Jacobi equation along a 
geodesic $\gamma$ in a semi-Riemannian manifold $(\M,\g)$
by means of a parallel trivialization of the tangent bundle $T\M$ along $\gamma$. 
Using such a trivialization, we may then identify vector fields along $\gamma$ with
curves in $\R^n$ and the metric tensor $\g$ along $\gamma$ with a fixed
nondegenerate symmetric bilinear form $g$ in $\R^n$. For all $t\in[a,b]$,
the endomorphism $v\mapsto \mathcal R(\dot\gamma(t),v)\,\dot\gamma(t)$  of
$T_{\gamma(t)}\M$ is identified with a $g$-symmetric endomorphism $R(t)$
of $\R^n$. Since covariant derivative along $\gamma$ corresponds to the
usual derivative of curves in $\R^n$, the Jacobi equation along
$\gamma$ becomes the Morse--Sturm system \eqref{eq:MS}.

If $\Pc$ is a smooth submanifold of $\M$ such that $\gamma(a)\in \Pc$
and $\dot\gamma(a)\in T_{\gamma(a)}\Pc^\perp$, then the tangent space
$T_{\gamma(a)}\Pc$ is identified with a subspace $P$ of $\R^n$, and the
second fundamental form $\mathcal S^\Pc_{\dot\gamma(a)}$ is identified
with a symmetric bilinear form $S$ on $P$. We assume that $\Pc$ is nondegenerate
at $\gamma(a)$, so that $g$ is nondegenerate on $P$. 

The space $\J$ of $\Pc$-Jacobi fields corresponds to the space $\Jb$ of 
solutions of \eqref{eq:MS} satisfying the initial conditions:
\begin{equation}\label{eq:ICMS}
v(a)\in P,\quad v'(a)+S\big(v(a)\big)\in P^\perp,
\end{equation}
where $\perp$ denotes the orthogonal complement with respect to $g$ and $S$ is
seen as a $g$-symmetric linear endomorphism of $P$.

We denote by $L^2([a,b];\R^m)$ the Hilbert space of square integrable 
$\R^m$-valued functions on $[a,b]$, by $H^1([a,b];\R^m)$ the Sobolev
space of absolutely continuous maps with derivative in $L^2([a,b];\R^m)$, and
by $H^1_0([a,b];\R^m)$ the subspace of $H^1([a,b];\R^m)$ consisting of functions
vanishing at $a$ and at $b$. We also denote by $C^0([a,b];\R^m)$ the 
Banach space of continuous functions from $[a,b]$ to $\R^m$.
It is well known that the inclusion maps $H^1([a,b];\R^m)\hookrightarrow C^0([a,b];\R^m)$
and $H^1([a,b];\R^m)\hookrightarrow L^2([a,b];\R^m)$ are compact operators
(see for instance \cite{Brezis}).

The Hilbert space $\mathcal H^\Pc_\gamma$ corresponds by the parallel trivialization
to the subspace $\mathcal H\subset H^1([a,b];\R^n)$ given by:
\begin{equation}\label{eq:defH}
\mathcal H=\Big\{v\in H^1([a,b];\R^n):v(a)\in P,\ v(b)=0\Big\};
\end{equation}
moreover, the index form $I^\Pc_\gamma$ defines a bounded symmetric bilinear
form $I$ on $\mathcal H$ by:
\begin{equation}\label{eq:defI}
I(v,w)=\int_a^b\Big[g(v',w')+g(Rv,w)\Big]\;\mathrm dt-S\big(v(a),w(a)\big).
\end{equation}
Observe that the kernel of $I$ in $\mathcal H$ is the space:
\begin{equation}\label{eq:kerI}
\mathrm{Ker}(I)=\mathcal H\cap\Jb.
\end{equation}

The notions of focal instants, multiplicity, signature, focal index and
Maslov index may be defined for Morse--Sturm systems \eqref{eq:MS}
with initial conditions \eqref{eq:ICMS} in the obvious way. 
\begin{defin}\label{thm:deffocMS}
An instant $t\in\left]a,b\right]$ is said to be {\em focal\/} for the Morse--Sturm
system \eqref{eq:MS} (with initial conditions \eqref{eq:ICMS}) if
there exists a non zero solution $v\in\Jb$ such that $v(t)=0$. 
The dimension of the space of such solutions is the {\em multiplicity\/}
of the focal instant. The {\em signature\/} of the focal instant $t$ is
defined to be the signature of the restriction of $g$ to $\Jb[t]^\perp$,
where:
\[\Jb[t]=\Big\{J(t):J\in\Jb\Big\}.\]
A focal instant is {\em nondegenerate\/} if $g$ is nondegenerate on $\Jb[t]$.
If there are only a finite number of focal instants, we define the
{\em focal index\/} of the Morse--Sturm system to be the sum
of the signatures of the focal instants in $]a,b]$. If $t=b$ is
not a focal instant, we define the {\em Maslov index\/} of the Morse--Sturm
system to be the number $\mu_{L_0}(\ell\vert_{[a+\varepsilon,b]})$, where
$\varepsilon>0$ is such that there are no focal instants in $]a,a+\varepsilon]$
and $\ell$, $L_0$ are defined in \eqref{eq:defellt} and \eqref{eq:L0}.
\end{defin}
Proposition~\ref{thm:imaslovifoc} generalizes in an obvious way to Morse--Sturm systems.

We are going
to prove a version of Theorem~\ref{thm:indexth} for such systems, which
in particular implies that the result holds in the geometrical
context. 

As a matter of facts, it is not hard to prove that every Morse--Sturm
system \eqref{eq:MS} with smooth coefficients arises from the Jacobi
equation along a semi-Riemannian geodesic, provided that one considers
a parallel trivialization of the normal bundle along the geodesic.
Details are found in \cite[Proposition~2.3.1]{MPT}.

Let us consider now a maximal negative distribution $\mathcal D$
along $\gamma$; each subspace $\mathcal D_t\subset T_{\gamma(t)}\M$ corresponds
to a subspace $D_t\subset\R^n$ by means of the parallel trivialization of $T\M$
along $\gamma$. Obviously, each $D_t$ is a maximal negative subspace for the
bilinear form $g$. 

The subspaces $\KgD$ and $\SgD $ of $\mathcal H^\Pc_\gamma$ correspond
to the closed subspaces $\Kt $ and $\St$ of $\mathcal H$ given by:
\begin{equation}\label{eq:KSt}
\begin{split}
\Kt=\Big\{v\in\mathcal H :&\;g(v',Y_i)\in
H^1([a,b];\R),\\&\;g(v',Y_i)'= g(v',Y_i')+g(Rv,Y_i),\ i=1,\ldots,k\Big\}
\\
\St=\Big\{v\in\mathcal H :&\;v(a)=0,\ v(t)\in D_t,\ \
\forall\,t\in[a,b]\Big\},
\end{split}
\end{equation}
where $Y_1,\ldots,Y_k$ is a frame for $D$, i.e., each $Y_i:[a,b]\to\R^n$ is a smooth
curve and $\{Y_1(t),\ldots,Y_k(t)\}$ is a basis of  $D_t$ for all $t$.

We are interested in determining the elements of the intersection $\Kt\cap\St$;
such elements are characterized as solutions of a second order linear
differential equation in $\R^n$ which is in general {\em not\/} a Morse--Sturm system.
This equation belongs to the more general class of {\em symplectic differential
systems}, that will be discussed in the next subsection.
\subsection{Symplectic differential systems in $\mathbf{I\!R^n}$}\label{sub:symplectic}
A Morse--Sturm system \eqref{eq:MS} can be written as the following first order linear
system in $\R^n\oplus{\R^n}^*$:
\begin{equation}\label{eq:MS1st}
\left( \begin{array}{c}v\\ \alpha \end{array}\right)'=
\left( \begin{array}{cc}0&g^{-1}\\ gR&0 \end{array}\right)
\left( \begin{array}{c}v\\ \alpha \end{array}\right),
\quad v(t)\in\R^n,\ \alpha(t)\in{\R^n}^*,
\end{equation}
where again the bilinear form $g$ is seen as a linear map from $\R^n$ to ${\R^n}^*$.

We denote by $\mathrm{Sp}(2n,\R)$ the Lie group of symplectic transformations
of the space $(\R^n\oplus{\R^n}^*,\omega)$, where $\omega$ is the symplectic form
defined in \eqref{eq:defomega}, and by
$\mathrm{sp}(2n,\R)$ its Lie algebra. 
Recall that an element $X\in\mathrm{sp}(2n,\R)$ is a linear endomorphism of 
$\R^n\oplus{\R^n}^*$ such that $\omega(X\,\cdot\,,\,\cdot\,)$
is symmetric; in block matrix form, $X$ is given by:
\begin{equation}\label{eq:X}
X=\left(\begin{array}{cc}A&B\\ C&-A^*\end{array}\right),
\end{equation}
where $A:\R^n\to\R^n$ is an arbitrary linear map, and $B:{\R^n}^*\to\R^n$,
$C:{\R^n}\to{\R^n}^*$ are symmetric when regarded as bilinear forms.

We observe that the coefficient matrix of the Morse--Sturm system
\eqref{eq:MS1st} is of the form \eqref{eq:X} with $A=0$, $B=g^{-1}$ and $C=gR$; 
we call a {\em symplectic  differential system\/} in $\R^n$ a first order linear
differential system in $\R^n\oplus{\R^n}^*$ whose coefficient matrix $X(t)$ is a
continuous curve in $\mathrm{sp}(2n,\R)$, where the blocks $A$ and $B$ are of class
$C^1$, and $B(t)$ is invertible for all $t\in[a,b]$:
\begin{equation}\label{eq:sympl}
\left\{\begin{array}{l}v'(t)=A(t)v(t)+B(t)\alpha(t);\\\alpha'(t)=C(t)v(t)-A^*(t)\alpha(t),
\end{array}\right.\quad t\in [a,b],\ v(t)\in\R^n,\ \alpha(t)\in{\R^n}^*.
\end{equation}

Morse--Sturm systems are special cases  of symplectic differential
systems with $A=0$ and $B$ constant; the index theory for Morse--Sturm
systems extends naturally to the class of symplectic differential
systems (see \cite{PT2}).  Such systems appear naturally as linearizations
of the Hamilton  equations, and also as the Jacobi equations along
geodesics when a non parallel trivialization of the tangent bundle is chosen.
Moreover, the class of symplectic differential systems is the more
natural class for which it is possible to define the notion of Maslov index
(see Section~\ref{sec:indexthsympl}).

We will need the extension of the index theory to symplectic systems
in order to calculate the term $n_+\left(I_\gamma^\Pc\big\vert_{\SgD}\right)$ 
that appears in equation  \eqref{eq:indexth}. Namely, the restriction
$I_\gamma^\Pc\big\vert_{\SgD}$ can be thought as the index form associated
to a symplectic system which is determined by the Jacobi equation along the
geodesic and by the choice of the distribution $\mathcal D$.

To clarify the situation, we outline briefly the basics of the index theory
for symplectic differential systems. 

Consider the symplectic differential system \eqref{eq:sympl} with coefficient
matrix $X$ given by \eqref{eq:X};  we say that a $C^1$-curve $v:[a,b]\to\R^n$
is an {\em $X$-solution\/} if there exists $\alpha:[a,b]\to{\R^n}^*$ of 
class $C^1$ such that the pair $(v,\alpha)$ is a solution of \eqref{eq:sympl}.
It is easy to see that an $X$-solution $v$
is of class $C^2$, and that, since $B$ is invertible, the unique $\alpha=\alpha_v$
such that $(v,\alpha)$ is a solution of \eqref{eq:sympl} is given by:
\begin{equation}\label{eq:alphav}
\alpha_v=B^{-1}(v'- Av).
\end{equation}
We denote by $\mathbb V$ the set of all $X$-solutions vanishing at $t=a$:
\begin{equation}\label{eq:spaces}
\mathbb V=\Big\{v: v\ \text{is an $X$-solution, with}\ v(a)=0\Big\}.
\end{equation}
Using the symmetry of $B$ and $C$ and \eqref{eq:sympl}, 
it is easy to see that the following equality holds:
\begin{equation}\label{eq:abvwcosnt2}
\alpha_v(w)=\alpha_w(v),\quad\forall\,v,w\in\mathbb V.
\end{equation}
For $t\in[a,b]$, we set 
\[\mathbb V[t]=\Big\{v(t):v\in\mathbb V\Big\}.\]
From \eqref{eq:abvwcosnt2} and a simple dimension counting argument,
the annihilator of $\mathbb V[t]$ is given by:
\begin{equation}\label{eq:VTo}
\mathbb V[t]^o=\Big\{\alpha_v(t):v\in\mathbb V,\ v(t)=0\Big\},\quad t\in\,[a,b]. 
\end{equation} 
\begin{defin}\label{thm:deffocal}
An instant $t\in\left[a,b\right]$ is said to be {\em focal\/} if there exists a non zero
$v\in\mathbb V$  such that $v(t)=0$, i.e., if $\mathbb V[t]\ne\R^n$.
The {\em multiplicity\/} $\mul(t)$ of the focal instant $t$ is defined to
be the dimension of the space of those $v\in\mathbb V$ vanishing at $t$, or,
equivalently, the codimension of $\mathbb V[t]$ in $\R^n$.
The {\em signature} $\sgn(t)$ of the focal instant $t$ is the signature
of the restriction of the bilinear form $B(t)$ to the space $\mathbb V[t]^o$, or,
equivalently, the signature of the restriction of $B(t)^{-1}$ to
the $B(t)^{-1}$-orthogonal complement $\mathbb V[t]^\perp$ of $\mathbb V[t]$ 
in $\R^n$.
The focal instant $t$ is said to be {\em nondegenerate\/} if such restriction
is nondegenerate.
If there is only a finite number of focal instants in $\left]a,b\right]$, we define
the {\em focal index\/} $\mathrm i_{\mathrm{foc}}=\mathrm i_{\mathrm{foc}}(X)$ 
to be the sum:
\begin{equation}\label{eq:defifocal}
\mathrm i_{\mathrm{foc}}=\sum_{t\in\left]a,b\right]}\sgn(t).
\end{equation}
\end{defin}
In the special situation that the bilinear form $B(t)$ is positive definite,
then the focal instants are obviously nondegenerate, and their
signatures   coincide  with their multiplicities. Moreover, as in
Proposition~\ref{thm:isolPfocalpts}, it can be proven that nondegenerate focal instants
are isolated.

The {\em index form\/} $I_X$ associated to \eqref{eq:sympl} is the bounded
symmetric bilinear form on the Hilbert space $H^1_0([a,b];\R^n)$
given by:
\begin{equation}\label{eq:IX}
I_X(v,w)=\int_a^b\Big[B(\alpha_v,\alpha_w)+C(v,w)\Big]\;\mathrm dt.
\end{equation}
Observe that, for a symplectic system \eqref{eq:MS1st} coming 
from a Morse--Sturm system \eqref{eq:MS}, the index form $I_X$ coincides
with the index form $I$ of formula \eqref{eq:defI} when the subspace
$P$ is chosen equal to $\{0\}$.

Integration by parts in \eqref{eq:IX} shows that the kernel of
$I_X$ is given by:
\begin{equation}\label{eq:kerIX}
\mathrm{Ker}(I_X)=\Big\{v\in\mathbb V:v(b)=0\Big\}.
\end{equation}
\begin{rem}\label{thm:remdefpos}
Observe that, from \eqref{eq:kerIX} we obtain easily that if
$B$ is positive definite and $C$ is positive semi-definite, then
the symplectic differential system \eqref{eq:sympl} has no focal instants.
\end{rem}
There is a natural notion of isomorphism in the class of symplectic
systems. 
Let $L_0$ be the Lagrangian subspace $\{0\}\oplus{\R^n}^*$ of $(\R^n\oplus{\R^n}^*,
\omega)$; we denote by $\mathrm{Sp}(2n,\R;L_0)$ the closed  subgroup of 
$\mathrm{Sp}(2n,\R)$ consisting of those symplectomorphisms $\phi_0$ such
that $\phi_0(L_0)=L_0$. It is easily seen that any such symplectomorphism
is given in block matrix form by:
\begin{equation}\label{eq:phi0}
\phi_0=\left(\begin{array}{cc}Z&0\\ { Z^*}^{-1}W&{
Z^*}^{-1}\end{array}\right),
\end{equation} with $Z:\R^n\to\R^n$ an isomorphism and $W$ a symmetric
bilinear form in $\R^n$.

We give the following:
\begin{defin}\label{thm:defequiv}
The symplectic differential systems with coefficient matrices $X$ and $\tilde X$
 are said to
be {\em isomorphic\/} if there exists a $C^1$-map $\phi_0:[a,b]\to
\mathrm{Sp}(2n,\R;L_0)$ whose upper-left $n\times n$ block is of class $C^2$ and
such that:
\begin{equation}\label{eq:isoX}
\tilde X=\phi_0'\phi_0^{-1}+\phi_0X\phi_0^{-1}.
\end{equation}
We call the map $\phi_0$ an {\em isomorphism\/} between $X$ and
$\tilde X$.
\end{defin}
The motivation of such notion of isomorphism is that, for
isomorphic systems $X$ and $\tilde X$, a pair $(v,\alpha)$
is an $X$-solution if and only if $\phi_0(v,\alpha)$ is an $\tilde X$-solution.
More precisely, we have the following relations between isomorphic
symplectic systems:
\begin{prop}\label{thm:isosympl}
Let $X$ and $\tilde X$ be the coefficient matrices of isomorphic symplectic systems,
and let $\phi_0$ as in formula \eqref{eq:phi0}
be an isomorphism between $X$ and $\tilde X$.

Then, the focal instants corresponding to the systems associated to $X$ and
$\tilde X$ are the same, and they have the same multiplicities and signatures.
Moreover, the isomorphism $v\mapsto Zv$ of $H^1_0([a,b];\R^n)$ carries
the index form $I_X$ into the index form $I_{\tilde X}$.
\end{prop}
\begin{proof}
See \cite[Subsection~2.10, Proposition~2.10.3]{PT2}.
\end{proof}
Although an index theory for symplectic systems may be developed
directly, the easiest way to extend the Morse Index theorem to this
class of systems is given by considering the following result:
\begin{prop}\label{thm:isoMSsympl}
Every symplectic system \eqref{eq:sympl} such that $B$ is a map of class
$C^2$ is isomorphic to a Morse--Sturm system \eqref{eq:MS1st}.
\end{prop}
\begin{proof}
See \cite[Proposition~2.11.1]{PT2}.
\end{proof}
Observe that the index $n_-(B)$ is invariant by isomorphisms of symplectic systems.
Hence we have the following:
\begin{cor}\label{thm:thindexsympl}
Consider the symplectic system \eqref{eq:sympl}, with $B$ a map of class
$C^2$ and {\em positive definite}. Then, there are only finitely many
focal instants, the index of $I_X$ in $H^1_0([a,b];\R^n)$ is finite, 
and it is equal to the sum of the multiplicities of the focal instants in $]a,b[$.
\end{cor}
\begin{proof}
The result is well known for Morse--Sturm systems 
(see for instance~\cite[Corollary~3.7]{GMPT}). The conclusion follows from 
Proposition~\ref{thm:isosympl} and Proposition~\ref{thm:isoMSsympl}.
\end{proof}

\subsection{The reduced symplectic system}\label{sub:reduced}
We now go back to the setup of Subsection~\ref{sub:reduction}
and we study the intersection of the spaces $\Kt$ and $\St$.

\begin{lem}\label{thm:systred}
Let $v\in\St$; write $v=\sum_{i=1}^kf_i\,Y_i$. Then, $v\in\Kt$ if and only
if $f=(f_1,\ldots,f_k)$ is a solution of the following symplectic
differential system:
\begin{equation}\label{eq:symplassoc}
\left\{\begin{array}{l}f'=-\mathcal B^{-1}\mathcal C\,f-\mathcal
B^{-1}\,\varphi,\\ \varphi'=(\mathcal C^*\mathcal B^{-1}\mathcal
C-\mathcal I)\,f
+\mathcal C^*\mathcal B^{-1}\,\varphi.\end{array}\right.
\end{equation}
where $\mathcal B$, $\mathcal I$
are bilinear forms in $\R^k$, and $\mathcal C$ 
is a linear map from $\R^k$ to ${\R^k}^*$, whose
matrices in the canonical basis are given by:
\begin{equation}
\label{eq:defcalobj}
\mathcal B_{ij}=g(Y_i,Y_j),\quad \mathcal C_{ij}=g(Y_j',Y_i),
\quad
\mathcal I_{ij}=g(Y_i',Y_j')+g(RY_i,Y_j).
\end{equation}
\end{lem}
\begin{proof}
It is a simple calculation based on the definition of the space $\Kt$ 
given in \eqref{eq:KSt}.
\end{proof}
\begin{cor}\label{thm:multfocinstred}
The dimension of the intersection $\Kt\cap\St$ is equal to the multiplicity
of $t=b$ as a focal instant for the symplectic differential system \eqref{thm:systred}.
\qed\end{cor}
\begin{defin}\label{thm:defsymplred}
The system \eqref{eq:symplassoc} is called the {\em reduced symplectic system\/}
associated to the Morse--Sturm system \eqref{eq:MS}, the maximal negative distribution
$D$ and the frame $Y_1,\ldots,Y_k$.
\end{defin}
It is not hard to prove that different choices of a frame for the distribution
$D$ produce isomorphic reduced symplectic systems (see \cite[Proposition~2.10.4]{PT2}).
\begin{rem}\label{thm:remisoredsympl}
It is easily seen that the following symplectic differential system
is isomorphic to \eqref{eq:symplassoc}:
\[
\left\{\begin{array}{l}f'=-\mathcal B^{-1}\mathcal C_{\mathrm a}f+\mathcal B^{-1}\varphi;\\
\varphi'=(\mathcal I-\mathcal C_{\mathrm s}'+\mathcal C_{\mathrm a}\mathcal B^{-1}
\mathcal C_{\mathrm a})f-\mathcal C_{\mathrm a}\mathcal B^{-1}\varphi,
\end{array}\right.\quad t\in [a,b],\ f(t)\in\R^k,\ \varphi(t)\in{\R^k}^*,
\]
where $\mathcal C_{\mathrm a}$ and $\mathcal C_{\mathrm s}$ are given by:
\begin{equation}\label{eq:defCaCs}
\mathcal C_{\mathrm a}=\frac12(\mathcal C-\mathcal C^*),\quad
\mathcal C_{\mathrm s}=\frac12(\mathcal C+\mathcal C^*).
\end{equation}
\end{rem}

The index form of the reduced symplectic system \eqref{eq:symplassoc}
corresponds to the restriction of $-I$ to the space $\St$:
\begin{prop}\label{thm:restrI}
The Hilbert space isomorphism \[H^1_0([a,b];\R^k)\ni
f=(f_1,\ldots,f_k)\longmapsto\sum_{i=1}^kf_i\cdot
Y_i\in\St\phantom{H^1_0([a,b];\R^k)\ni}\] carries the index form of the reduced
symplectic system \eqref{eq:symplassoc} to the restriction of $-I$ to $\St$, where $I$
is the index form of the original Morse--Sturm system defined in \eqref{eq:defI}.
\end{prop}
\begin{proof}
It is an easy calculation that uses \eqref{eq:IX}.
\end{proof}
\begin{cor}\label{thm:indIS}
The coindex $n_+(I\vert_{\St})$ of the restriction of $I$ to
$\St$ is finite, and it is equal to the sum of the multiplicities of the
conjugate instants of the reduced symplectic system \eqref{eq:symplassoc}
in $]a,b[$.
\end{cor}
\begin{proof}
Observe that the coefficient of $\varphi$ in the first equation of
\eqref{eq:symplassoc} is positive definite. The conclusion follows
from Corollary~\ref{thm:thindexsympl} and Proposition~\ref{thm:restrI}.
\end{proof}
We now give a criterion for the vanishing of the number $n_+(I\vert_{\St})$:
\begin{cor}\label{thm:corcriterio}
Suppose that either one of the following symmetric bilinear
forms (see \eqref{eq:defcalobj} and \eqref{eq:defCaCs}):
\[\mathcal C^*\mathcal B^{-1}\mathcal C-\mathcal I,\quad \mathcal C_{\mathrm s}'-
\mathcal C_{\mathrm a}\mathcal B^{-1}\mathcal C_{\mathrm a}-\mathcal I\]
is positive semi-definite on $[a,b]$. Then $n_+(I\vert_{\St})=0$.
\end{cor}
\begin{proof}
It follows directly from Remark~\ref{thm:remdefpos}, Remark~\ref{thm:remisoredsympl}
and Corollary~\ref{thm:indIS}. 
\end{proof}

\begin{cor}\label{thm:Idefneg}
The restriction of $I$ to $\St$ is represented by a self-adjoint
operator on $\St$ which is a compact perturbation of a negative
isomorphism of $\St$.
\end{cor}
\begin{proof}
The index form of any symplectic differential system \eqref{eq:sympl}
with the coefficient $B$ positive definite is represented
by a compact perturbation of a positive isomorphism of 
 $H^1_0([a,b];\R^n)$ (see \cite[Lemma~2.6.6]{PT2}). The conclusion follows from
Proposition~\ref{thm:restrI}.
\end{proof}
\subsection{An extension of the index form}\label{sub:sustenido}
The strategy for proving Theorem~\ref{thm:indexth} will be to apply
Proposition~\ref{thm:notsoelementary} to a family $I_t$ of symmetric bilinear
forms on Hilbert spaces  $\Kt_t$  obtained by considering restrictions
of the Morse--Sturm system \eqref{eq:MS} to the interval $[a,t]$.
Unfortunately, we run into the annoying technical problem that the family
$\Kt_t$ fails to be $C^1$ around the focal instants of the reduced
symplectic system \eqref{eq:symplassoc}. 

In this subsection we describe a trick to overcome this problem by
introducing an artificial extension $I^\#$ of the index form $I$
to a space $\Kt^\#$ so that the corresponding family $\Kt^\#_t$
will be of class $C^1$.

Let us introduce the ``sharped'' versions of the objects of our theory:
$\mathcal H^\#$, $\Kt^\#$, $\St^\#$ and $I^\#$.
Set:
\begin{equation}\label{eq:sustobjects}
\begin{split}
\mathcal H^\#=&\;\Big\{v\in H^1([a,b];\R^n):v(a)\in P\Big\},\\
\Kt^\#=&\;\Big\{v\in\mathcal H^\# : g(v',Y_i)\in
H^1([a,b];\R),\\&\phantom{\;\Big\{v\in\mathcal H^\# :}\;g(v',Y_i)'=
g(v',Y_i')+g(Rv,Y_i),\ i=1,\ldots,k\Big\},
\\
\St^\#=&\;\Big\{v\in\mathcal H^\# : v(a)=0,\ v(t)\in D_t,\ \
\forall\,t\in[a,b]\Big\}.
\end{split}
\end{equation}
 Throughout this subsection
we consider a fixed symmetric bilinear form $\Theta$ in $\R^n$.
The extended bilinear form $I^\#$ is defined using $\Theta$ by:
\begin{equation}\label{eq:sustI}
I^\#(v,w)=\int_a^b\Big[g(v',w')+g(Rv,w)\Big]\;\mathrm
dt+\Theta\big(v(b),w(b)\big)- S\big(v(a),w(a)\big).
\end{equation}
Recall that $\Jb$ denotes the space of solutions of \eqref{eq:MS}
satisfying the initial conditions \eqref{eq:ICMS}; we can characterize
the kernel of $I^\#$ as:
\begin{equation}\label{eq:KerIsust}
\mathrm{Ker}(I^\#)=\Big\{v\in\Jb:g\,v'(b)+\Theta\big(v(b)\big)=0\Big\},
\end{equation}
where $g$ and $\Theta$ are considered as linear maps from $\R^n$ to ${\R^n}^*$.

Observe that $\Kt=\Kt^\#\cap\mathcal H$, $\St=\St^\#\cap \mathcal H$ and
that $I$ is the restriction of $I^\#$ to $\mathcal H$.
\smallskip

We define a bounded linear map $F:\mathcal H^\#\to L^2([a,b];{\R^k}^*)$:
\begin{equation}\label{eq:defF}
\big[F(v)
(t)\big]_i=g\big(v'(t),Y_i(t)\big)-\int_a^t\Big[g(v',Y_i')+g(Rv,Y_i)\Big]\;\mathrm ds,
\quad i=1,\ldots,k.
\end{equation}
Obviously, $\Kt^\#$ is the inverse image by $F$ of the subspace $\mathfrak C $ of
$L^2([a,b];{\R^k}^*)$ consisting of constant functions.

\begin{lem}\label{thm:Fiso}
The restriction of $F$ to $\St^\#$ is an isomorphism.
\end{lem}
\begin{proof}
We identify the space $\St^\#$ with the  space $\mathcal X=\{f\in
H^1([a,b];\R^k):f(a)=0\}$ by  the map $v=\sum_{i}f_i\,Y_i\mapsto f=(f_1,\ldots,f_k)$;
then using
\eqref{eq:defcalobj}, the map $F$ on $\St^\#$ can be written as:
\begin{equation}\label{eq:Fwithcal}
F(v)(t)=\mathcal B(f')(t)+\mathcal C(f)(t)-\int_a^t\big[\mathcal
C^*(f')+\mathcal I(f)\big]\mathrm ds.
\end{equation}
Using the fact that the inclusion of $H^1$ in $L^2$ is compact, it is easy to see
from \eqref{eq:Fwithcal} that the restriction of $F$ to
$\St^\#$ is a compact perturbation of the isomorphism
$\mathcal X\ni f\mapsto \mathcal B(f')\in L^2([a,b];{\R^k}^*)$.
Hence, the restriction of $F$ to $\St^\#$ is a Fredholm operator of index zero,
and to prove the Lemma it suffices to show that $F$ is injective
on $\St^\#$. To this aim, observe that if $f\in\mathrm{Ker}(F\vert_{\St^\#})$,
then using \eqref{eq:Fwithcal} we see that $f$ is a solution of a second
order homogeneous linear differential equation, and that $f(a)=f'(a)=0$.
This concludes the proof.
\end{proof}
\begin{cor}\label{thm:Fsobre}
The map $F$ is surjective, $\mathrm{dim}(\Kt^\#\cap\St^\#)=k$ and
$\mathcal H^\#=\Kt^\#+\St^\#$.
\end{cor} 
\begin{proof}
It follows easily from Lemma~\ref{thm:Fiso} and the fact that $\Kt^\#$ is the inverse
image by $F$ of the space $\mathfrak C$ 
of constant functions, which is $k$-dimensional.
\end{proof}
\begin{lem}\label{thm:Fsobreconst}
Let $q: L^2([a,b];{\R^k}^*)\to L^2([a,b];{\R^k}^*)/\mathfrak C$ be the quotient
map. Suppose that $\Kt\cap\St=\{0\}$; then,   $q\circ F$ maps $\St$
isomorphically onto $  L^2([a,b];{\R^k}^*)/\mathfrak C$.
\end{lem}
\begin{proof}
The proof is essentially the same as the proof of Lemma~\ref{thm:Fiso}.
Namely, using \eqref{eq:Fwithcal} we show that the restriction of $q\circ F$
to $\St$ is a Fredholm operator of index zero. The injectivity of this
restriction is obviously equivalent to $\Kt\cap\St=\{0\}$.
\end{proof}
\begin{cor}\label{thm:qFsurjective}
If $\Kt\cap\St=\{0\}$, then $q\circ F$ is surjective and $\mathcal H=\Kt\oplus\St$. 
\end{cor}
\begin{proof}
It follows from Lemma~\ref{thm:Fsobreconst} and the fact that
$\Kt=\mathrm{Ker}\left(q\circ F\vert_{\mathcal H}\right)$.
\end{proof}
\begin{lem}\label{thm:codim}
The space $\left(\Kt^\#+\St\right)$ is closed in
$\mathcal H^\#$, and its codimension is equal to $\mathrm{dim}(\Kt \cap\St)$.
\end{lem}
\begin{proof}
The fact that $\left(\Kt^\#+\St\right)$ is closed follows easily from
$\mathcal H^\#=\Kt^\#+\St^\#$ (Corollary~\ref{thm:Fsobre}).
Now, we observe:
\begin{equation}\label{eq:codim1}
\mathrm{codim}\left(\Kt^\#+\St\right)=\mathrm{dim}\left(\frac{\mathcal
H^\#}{\Kt^\#+\St}\right).
\end{equation}
We have a surjective map:
\begin{equation}\label{eq:codim2}
\iota:\frac{\St^\#}\St\longrightarrow\frac{\Kt^\#+\St^\#}{\Kt^\#+\St}=\frac{\mathcal
H^\#}{%
\Kt^\#+\St}
\end{equation}
induced by inclusion, and 
\begin{equation}\label{eq:codim3}
\mathrm{Ker}(\iota)=\frac{\St^\#\cap(\Kt^\#+\St)}{\St}=\frac{(\St^\#\cap\Kt^\#)+\St}{%
\St}\simeq\frac{\St^\#\cap\Kt^\#}{\Kt\cap\St}.
\end{equation}
Finally, using Corollary~\ref{thm:Fsobre}, \eqref{eq:codim1}, \eqref{eq:codim2}
and \eqref{eq:codim3} we compute:
\begin{equation}\label{eq:codim4}
\begin{split}
\mathrm{codim}(\Kt^\#+\St)&\;=\mathrm{dim}(\mathrm{Im}(\iota))=\mathrm{dim}\left(
\frac{\St^\#}\St\right)-\mathrm{dim}\left(\frac{\St^\#\cap\Kt^\#}{\Kt\cap\St
}\right)=\\&= k-\big(k-\mathrm{dim}(\Kt \cap\St)\big)=\mathrm{dim}(\Kt \cap\St).
\qedhere\end{split}
\end{equation}
\end{proof}
\begin{lem}\label{thm:KSorth}
The spaces $\Kt^\#$ and $\St$ are $I^\#$-orthogonal, i.e., $I^\#(v,w)=0$ for
all $v\in\Kt^\#$ and all $w\in\St$.
\end{lem}
\begin{proof}
Let $v\in\Kt^\#$ and $w=\sum_if_iY_i\in\St$ be given, with $f_i\in H^1_0([a,b];\R)$.
Using the definition of $\Kt^\#$, we compute:
\[I^\#(v,w)=\sum_{i=1}^k\int_a^b \frac{\mathrm d}{\mathrm dt}\,\big[f_i\,g(
v',Y_i)\big] \;\mathrm dt=0.\qedhere\]
\end{proof}
\begin{prop}\label{thm:KerIK} If $\Kt\cap\St=\{0\}$, then
the kernel of the restriction of $I$ to $\Kt$ is equal to
the kernel of $I$ in $\mathcal H$, given in formula \eqref{eq:kerI}.
\end{prop}
\begin{proof}
From Lemma~\ref{thm:KSorth} it follows that $\Kt$ and $\St$ are $I$-orthogonal;
the conclusion follows from Corollary~\ref{thm:qFsurjective} and the observation that
the kernel of $I$ in $\mathcal H$ is contained in $\Kt$.
\end{proof}
\begin{prop}\label{thm:kerIsustK}
Suppose that $I^\#$ is nondegenerate on $\mathcal H^\#$. Then, the kernel of the
restriction of $I^\#$ to $\Kt^\#$ is given by $\Kt\cap\St$.
\end{prop}
\begin{proof}
The nondegeneracy assumption of $I^\#$ on $\mathcal H^\#$ means that it is
represented by an injective operator on $\mathcal H^\#$. Using the compact inclusion
of $H^1$ in $C^0$, it is easily seen that formula \eqref{eq:sustI} defines a bilinear
form which is represented by a compact perturbation of an isomorphism of
$H^1([a,b],\R^n)$. Using the additivity of the Fredholm index of operators it is easily
proven that $I^\#$ is represented by a Fredholm operator of index zero in $\mathcal
H^\#$; hence it follows that $I^\#$ is indeed represented by an isomorphism of $\mathcal
H^\#$.

Using Lemma~\ref{thm:KSorth}, we have inclusions:
\begin{equation}\label{eq:kerIK1}
\Kt\cap\St\subset\mathrm{Ker}\left(I^\#\big\vert_{\Kt^\#}\right)\subset
\Big\{v\in\mathcal H^\#:I^\#\big(v,\Kt^\#+\St\big)=0\Big\}.
\end{equation}
Since $I^\#$ is an isomorphism, the dimension of the third member
in \eqref{eq:kerIK1} equals the dimension of the annihilator
of $\Kt^\#+\St$ in $(\mathcal H^\#)^*$. The dimension of this annihilator
coincides with the codimension of $\Kt^\#+\St$ in $\mathcal H^\#$; by
Lemma~\ref{thm:codim}, it follows that the inclusions in \eqref{eq:kerIK1}
are equalities, which concludes the proof.
\end{proof}
\begin{prop}\label{thm:IdComp}
The restriction of $I^\#$ to $\Kt^\#$  (respectively, of $I$ to $\Kt$) 
is represented by a compact perturbation of a positive isomorphism of $\Kt^\#$
(respectively, of $\Kt$).
\end{prop}
\begin{proof}
It is essentially identical to the proof of \cite[Lemma~2.6.6]{PT2}.
\end{proof}
\begin{prop}\label{thm:indKsustindK}
Suppose that $t=b$ is not a focal instant for the Morse--Sturm system
\eqref{eq:MS}. Then,
\begin{equation}\label{eq:indKsustindK}
n_-\big(I^\#\big\vert_{\Kt^\#}\big)=n_-\Big(I\big\vert_{\Kt}\Big)+n_-\Big(\Theta-
\phi_{L_1,L_0}(\ell(b))\Big),
\end{equation}
where $L_0=\{0\}\oplus{\R^n}^*$, $L_1=\R^n\oplus\{0\}$, $\phi_{L_1,L_0}$
is the chart of the Lagrangian Grassmannian $\Lambda$ defined in \eqref{eq:defphiLoL1}
and $\ell:[a,b]\to\Lambda$ is the curve of Lagrangians defined in \eqref{eq:defellt}.
\end{prop}
\begin{proof}
Let us denote by $\beta$ the symmetric bilinear form on $L_1$ given
by $\phi_{L_1,L_0}(\ell(b))$; we regard $\beta$ as a linear map
from $\R^n$ to ${\R^n}^*$ by identifying $L_1\simeq\R^n$. By the definition
of the chart $\phi_{L_1,L_0}$, we have:
\[\ell(b)=\Big\{(v,\alpha):\alpha+\beta(v)=0\Big\};\]
therefore, we obtain:
\begin{equation}\label{eq:indKsustindK1}
\beta\big(v(b)\big)=-g\,v'(b),\quad\forall\,v\in\Jb.
\end{equation}
It is an easy observation that, since $t=b$ is not a focal
instant, we have:
\[\Kt^\#=\Kt\oplus\Jb,\]
where the direct sum is $I^\#$-orthogonal, and so:
\begin{equation}\label{eq:indKsustindK2}
n_-\big(I^\#\big\vert_{\Kt^\#}\big)=n_-
\Big(I\big\vert_{\Kt}\Big)+n_-\Big(I^\#\big\vert_\Jb\Big).
\end{equation}
The conclusion follows from the fact that, using \eqref{eq:indKsustindK1},
it is easily seen that the isomorphism $\Jb\ni v\mapsto v(b)\in\R^n$
carries the restriction of $I^\#$ to the bilinear form $\Theta-\beta$.
\end{proof}

\subsection{The index function $\mathbf{i(t)}$}
\label{sub:indexfunction}
In this subsection we consider the restriction of the Morse--Sturm
system \eqref{eq:MS} to the interval $[a,t]$, with $t\in\left]a,b\right]$.
We define the objects $\mathcal H_t$, $\mathcal H^\#_t$, $I_t$, $I_t^\#$,
$\Kt_t$, $\Kt^\#_t$, $\St_t$, $\St_t^\#$, $F_t$ as in formulas
\eqref{eq:defH}, \eqref{eq:defI}, \eqref{eq:KSt},
\eqref{eq:sustobjects}, \eqref{eq:sustI} and \eqref{eq:defF}  by replacing
$b$ with $t$. 
The definition of the bilinear form $I^\#_t$ depends on the choice
of a symmetric bilinear form $\Theta$; such choice will be made appropriately
when needed.

Clearly, all the results of the previous subsections remain valid
when the Morse--Sturm system is restricted to the interval $[a,t]$.
\smallskip

We study the evolution of the {\em index function}:
\begin{equation}\label{eq:defi}
i(t)=n_-\left(I_t\big\vert_{\Kt_t}\right),\quad t\in\left]a,b\right];
\end{equation}
obviously,
\[n_-\left(I\big\vert_{\Kt}\right)=i(b).\]
We will use the isomorphisms $\Phi_t:\mathcal H^\#\to\mathcal H_t^\#$ defined by:
$\Phi_t (\hat v)=v$, where 
\begin{equation}\label{eq:Phit}
v(s)=\hat v(u_s),\quad u_s=a+\frac{b-a}{t-a}(s-a),\quad\forall\,s\in[a,t];
\end{equation}
observe that $\Phi_t$ carries $\mathcal H$ onto $\mathcal H_t$.

We get families of closed subspaces of $\mathcal H^\#$ given by:
\[ 
\hat\Kt_t=\Phi_t^{-1}(\Kt_t),\quad \hat\Kt_t^\#=\Phi_t^{-1}\left(\Kt_t^\#\right),\quad
\hat\St_t=\Phi_t^{-1}(\St_t),\quad \hat\St_t^\#=\Phi_t^{-1}\left(\St_t^\#\right),
\]
we also get curves $\hat I:\left]a,b\right]\to\Bsym(\mathcal H,\R)$, $\hat
I^\#:\left]a,b\right]\to\Bsym(\mathcal H^\#,\R)$ of symmetric bilinear forms and 
a curve $\hat F:\left]a,b\right]\to \mathcal L\big(\mathcal H^\#,
L^2([a,b];{\R^k}^*)\big)$ of maps, defined by:
\[\hat I_t=I(\Phi_t\,\cdot\,,\Phi_t\,\cdot\,),\quad\hat
I_t^\#=I^\#(\Phi_t\,\cdot\,,\Phi_t\,\cdot\,), \quad\hat F_t=\Phi_t^{-1}\circ F_t\circ
\Phi_t.\]
We are also denoting by $\Phi_t$ the isomorphism from $L^2([a,b];{\R^k}^*)$
to $L^2([a,t];{\R^k}^*)$ defined by formula \eqref{eq:Phit}.

An explicit formula for $\hat I^\#_t$ is given by:
\begin{equation}\label{eq:explicit}
\begin{split}
\hat I^\#_t(\hat v,\hat w)=\;&\int_a^t\left[\left(\frac{b-a}{t-a}\right)^2
g\big(\hat v'(u_s),\hat w'(u_s)\big)+g\big(R(s)\hat v(u_s),\hat
w(u_s)\big)\right]\;\mathrm ds\\& +\Theta\big(\hat v(b),\hat w(b)\big)-S\big(\hat
v(a),\hat w(a)\big),
\end{split}
\end{equation}
for all $\hat v, \hat w\in\mathcal H^\#$.
\smallskip

As to the {\em initial value\/} $i(a)$ of the index function, we need
to consider suitable {\em extensions\/} to $t=a$ of the objects $\hat I_t$,
$\hat\Kt_t$ and $\hat F_t$. We set:
\[\mathfrak I_t=(t-a)  \hat I_t,\quad \mathcal F_t=(t-a)  \hat F_t,
\quad t\in\left]a,b\right].\]
A change of variable in \eqref{eq:explicit} gives the following explicit
formula for $\mathfrak I_t$:
\begin{equation}\label{eq:explicit2}
\begin{split}
\mathfrak I_t(\hat v,\hat w)=&\;\int_a^b\left[(b-a)\,g\big(\hat v'(u),\hat
w'(u)\big)+\frac{(t-a)^2}{b-a}\,g\big(R(s_u)\hat v(u),\hat w(u)\big)\right]\;\mathrm du\\
&+(t-a)\Big(\Theta\big(\hat v(b),\hat w(b)\big)-S\big(\hat v(a),\hat w(a)\big)\Big),
\end{split}
\end{equation}
for all $\hat v,\hat w\in\mathcal H$, where $s_u=a+\frac{t-a}{b-a}(u-a)$.
Setting $t=a$ in \eqref{eq:explicit2}, we define  $\mathfrak I_a$ as:
\begin{equation}\label{eq:Ia}
\mathfrak I_a(\hat v,\hat w)=(b-a)\int_a^b g\big(\hat v'(u),\hat
w'(u)\big)\;\mathrm du.
\end{equation}
Observe that:
\[\hat\Kt_t=\mathrm{Ker}\big(q\circ\hat F_t\big\vert_{\mathcal H}\big)=
\mathrm{Ker}\big(q\circ {\mathcal F}_t\big\vert_{\mathcal
H}\big),\quad\forall\,t\in\left]a,b\right],\] where $q:L^2([a,t];{\R^k}^*)\to
L^2([a,t];{\R^k}^*)/\mathfrak C$ is the quotient map and $\mathfrak C$ is the space of
constant functions.

An explicit formula for $\mathcal F_t$ is given by:
\begin{equation}\label{eq:explicit3}
\begin{split}
\Big[\mathcal F_t(\hat v)(u)\Big]_i=&\;(b-a)\,g\big(\hat v'(u),Y_i(s_u)\big)\\
&-\int_a^u\Big[(t-a)g\big(\hat v'(x),Y_i'(r_x)\big)+\frac{(t-a)^2}{b-a}\,
g\big(R(r_x)\hat v(x),Y_i(r_x)\big)\Big]\mathrm dx,
\end{split}
\end{equation}
$i=1,\ldots,k$, for all $\hat v\in\mathcal H^\#$, where
$r_x=a+\frac{s_u-a}{u-a}(x-a)$. Setting $t=a$ in \eqref{eq:explicit3}
gives the following definition for $\mathcal F_a$:
\begin{equation}\label{eq:defFa}
\Big[\mathcal F_a(\hat v)(u)\Big]_i=(b-a) g\big(\hat v'(u),Y_i(a)\big),\quad
\forall\,\hat v\in \mathcal H^\#.
\end{equation}
We also set $\hat\Kt_a=\mathrm{Ker}\big(q\circ\mathcal F_a\vert_{\mathcal H}\big)$,
namely,
\begin{equation}\label{eq:defKa}
\hat\Kt_a=\Big\{\hat v\in\mathcal H:g\big(\hat v'(u),Y_i(a)\big)=\text{constant},\ \ 
i=1,\ldots,k\Big\}.
\end{equation}

\begin{prop}\label{thm:tudoC1}
Suppose that $R$ is a map of class $C^1$. Then, $\hat I^\#:\left]a,b\right]\to
\Bsym(\mathcal H^\#,\R)$ and $\mathfrak I:[a,b]\to\Bsym(\mathcal H,\R)$ 
are maps of class $C^1$.  Moreover, $\big\{\hat\Kt^\#_t\big\}_{t\in\left]a,b\right]}$ is
a 
$C^1$-family of closed subspaces of $\mathcal H^\#$ and, provided that
there are no focal instants for the reduced symplectic system \eqref{eq:symplassoc} 
in the interval
$[c,d]\subset[a,b]$,  $\{\hat\Kt_t\}_{t\in[c,d]}$ is a $C^1$-family of closed 
subspaces of $\mathcal H$. 
\end{prop}
\begin{proof}
By standard regularity arguments (see \cite[Lemma~2.3, Proposition~3.3 and
Lem\-ma~4.3]{GMPT}), formula \eqref{eq:explicit2} shows that $\mathfrak I$ 
and $I^\#$ are $C^1$
in $[a,b]$, which obviously implies that $\hat I$ is $C^1$ in $]a,b]$.
Similarly, formula \eqref{eq:explicit3} shows that $\mathcal F$ is of
class $C^1$ on $[a,b]$; from Corollary~\ref{thm:Fsobre} we deduce that
$\mathcal F_t$ is surjective for $t\in\left]a,b\right]$.  
The regularity of the family $\big\{\hat\Kt^\#_t\big\}_{t
\in\left]a,b\right]}$ follows then from Lemma~\ref{thm:produce}. 

As to the regularity of the family $\{\hat\Kt_t\}_{t\in[c,d]}$, we have to show that
$q\circ \mathcal F_t\vert_{\mathcal H}$ is surjective for $t\in[c,d]$.
For $t=a$ it follows directly from the definition of $\mathcal F_a$
in \eqref{eq:defFa}. For $t>a$, the surjectivity follows from
Corollary~\ref{thm:multfocinstred} and
Corollary~\ref{thm:qFsurjective}.
\end{proof}
\begin{cor}\label{thm:iconstant}
Suppose that $R$ is a map of class $C^1$. If there are no focal instants
of the Morse--Sturm system \eqref{eq:MS} and also of the reduced
symplectic system \eqref{eq:symplassoc} in the interval $[c,d]\subset\, ]a,b]$,
then the index function $i$ is constant on $[c,d]$.
\end{cor}
\begin{proof}
Let $t\in[c,d]$ be fixed. Using formula
\eqref{eq:kerI}, Corollary~\ref{thm:multfocinstred} and Proposition~\ref{thm:KerIK}
we conclude that $I_t$ is nondegenerate on $\Kt_t$, and therefore $\hat I_t$
is nondegenerate on $\hat\Kt_t$. 
Keeping in mind the result of Proposition~\ref{thm:IdComp} and
Proposition~\ref{thm:tudoC1}, the conclusion follows  by applying 
Remark~\ref{thm:remcasoparticular}.
\end{proof}
\subsection{Proof of Theorem~\ref{thm:indexth}}\label{sub:proof}
We start with the following:
\begin{lem}\label{thm:contachata}
Let $t_0\in\left]a,b\right]$  and  $\hat v_0,\hat w_0\in\hat\Kt_{t_0}\cap\hat\St_{t_0}$
be fixed. Let $\hat{\mathfrak v},\hat{\mathfrak w}:\left]a,b\right]\to\mathcal H^\#$
be $C^1$-curves with $\hat{\mathfrak v}_t ,\hat{\mathfrak w}_t\in \hat\Kt^\#_t$
for all $t\in\left]a,b\right]$ and with $\hat{\mathfrak v}_{t_0}=\hat v_0$,
$\hat{\mathfrak w}_{t_0}= \hat w_0$. Then:
\begin{equation}\label{eq:contachata}
\frac{\mathrm d}{\mathrm dt}\,\hat I^\#_t(\hat{\mathfrak v}_t,\hat{\mathfrak w}_t)
\,\Big\vert_{t=t_0}=g\big(v_0'(t_0),w_0'(t_0)\big),
\end{equation}
where $v_0=\Phi_{t_0}(\hat v_0)$ and $w_0=\Phi_{t_0}(\hat w_0)$.
\end{lem}
\begin{proof}
By Remark~\ref{thm:remindepchoice}, the term on the left hand
side of \eqref{eq:contachata} does not depend on the choice of $\hat{\mathfrak v}$
and $\hat{\mathfrak w}$. In order to facilitate the computation, we make
a suitable choice of the curves $\mathfrak v_t$ and $\mathfrak w_t$, as follows.
Write $v_0(s)=\sum_{i=1}^kf_i^{(1)}(s)Y_i(s)$ and
$w_0(s)=\sum_{i=1}^kf_i^{(2)}(s)Y_i(s)$,
with $s\in[a,t_0]$; by Lemma~\ref{thm:systred}, the maps $\{f_i^{(1)}\}_i$ and
$\{f_i^{(1)}\}_i$ are solutions of the reduced symplectic system \eqref{eq:symplassoc},  
hence they define maps of class $C^2$ on the entire interval $[a,b]$.
We set \[\mathfrak v_t(s)=\sum_{i=1}^kf_i^{(1)}(s)Y_i(s),
\quad \mathfrak w_t(s)=\sum_{i=1}^kf_i^{(2)}(s)Y_i(s),\qquad s\in[0,t],\
t\in\left]a,b\right];\] again by Lemma~\ref{thm:systred}, $\mathfrak v_t$ and $\mathfrak
w_t$ are in $\mathcal K^\#_t$ for all $t$, and so the maps $\hat{\mathfrak v}_t$ and
$\hat{\mathfrak w}_t$ defined by:
\[\hat{\mathfrak v}_t=\Phi_t^{-1}(\mathfrak v_t),\qquad \hat{\mathfrak w}_t=
\Phi_t^{-1}(\mathfrak w_t)\]
are in $\hat{\mathcal K}_t^\#$. Obviously, the maps $(t,u)\mapsto\hat{\mathfrak v}_t(u)$
and $(t,u)\mapsto\hat{\mathfrak w}_t(u)$ are of class $C^2$, and therefore they define
$\mathcal H^\#$-valued $C^1$-maps.

Once the choice of $\hat{\mathfrak v}_t$ and $\hat{\mathfrak w}_t$ is
made, we  compute  as follows:
\begin{equation}\label{eq:chata1}
\begin{split}
\frac{\mathrm d}{\mathrm dt}\,\hat I^\#_t(\hat{\mathfrak v}_t,\hat{\mathfrak w}_t)
\,\big\vert_{t=t_0}= \frac{\mathrm d}{\mathrm dt}\, I^\#_t({\mathfrak
v}_t,{\mathfrak w}_t)\,\big\vert_{t=t_0}=g\big(v_0'(t_0),w_0'(t_0)\big),
\end{split}
\end{equation}
using that $v_0(t_0)=w_0(t_0)=0$ and that $\mathfrak v_t(s)$, $\mathfrak w_t(s)$
do not depend on $t$.
\end{proof}

\begin{cor}\label{thm:cornegdef}
Let $t_0\in\left]a,b\right]$ and suppose that $I^\#_{t_0}$ is nondegenerate
in $\mathcal H^\#_{t_0}$.
Setting $B(t)=\hat I^\#_t$ and $\mathcal D_t=\hat\Kt_t^\#$, then the
symmetric bilinear form $\overline B'(t_0)$ on $\mathrm{Ker}({\overline B}(t_0))=
\hat\Kt_{t_0}\cap\hat\St_{t_0}$ introduced in Definition~\ref{thm:defderivada} is
negative definite.
\end{cor}
\begin{proof}
Recall  that the kernel of  ${\overline B}(t_0)$ is given
in Proposition~\ref{thm:kerIsustK}. For $v_0,w_0\in\St_{t_0}$, 
we have $v_0'(t_0),w_0'(t_0)\in D_{t_0}$, and $g$ is negative definite
in $D_{t_0}$. From Lemma~\ref{thm:contachata} it follows that 
$\overline{B}'(t_0)$ is negative semi-definite. To conclude the
proof we have to show that, if $v_0\in\Kt_{t_0}\cap\St_{t_0}$ and
$v_0'(t_0)=0$, then $v_0=0$. This follows easily from Lemma~\ref{thm:systred}.
\end{proof}

We now determine the initial value of the index function $i(t)$.

\begin{lem}\label{thm:indJa}
The restriction of the symmetric bilinear form $\mathfrak I_a$ to $\hat{\Kt}_a$
is represented by a compact perturbation of a positive isomorphism. Moreover, it
is nondegenerate, and its index equals the index of the restriction of $g$
to $P$.
\end{lem}
\begin{proof}
See \cite[Lemma~2.7.8]{PT2}.
\end{proof}
\begin{cor}\label{thm:i(a)}
For $t\in\left]a,b\right]$ sufficiently close to $a$, we have $i(t)=n_-(g\vert_P)$.
\end{cor}
\begin{proof}
Obviously, for $t\in\left]a,b\right]$, $i(t)=n_-\big(\mathfrak
I_t\vert_{\hat\Kt_t}\big)$. The conclusion follows from
Remark~\ref{thm:remcasoparticular} and Lemma~\ref{thm:indJa}.
\end{proof}
We are finally ready for:
\begin{proof}[Proof of Theorem~\ref{thm:indexth}]
The proof will be done for Morse--Sturm systems (see Subsection~\ref{sub:reduction})
with coefficients of class $C^1$;
the geometrical version of the theorem is an immediate corollary.

We first consider the case that there are only a finite number 
of focal instants for the Morse--Sturm system \eqref{eq:MS} 
and that $t=b$ is not a focal instant for the reduced symplectic
system \eqref{eq:symplassoc}. Observe that, by Corollary~\ref{thm:thindexsympl},
the number of focal instants for the reduced symplectic system is finite.

By Corollary~\ref{thm:iconstant}, the function $i(t)$ is piecewise constant, with jumps
at the focal instants of either the Morse--Sturm system or the reduced
symplectic system. By Corollary~\ref{thm:i(a)}, $i(t)=n_-(g\vert_P)$
for $t$ sufficiently close to $a$. 

Let $t_0\in\left]a,b\right[$ be a focal instant for either the Morse--Sturm
or the reduced symplectic system; we compute the jump
of $i$ at $t_0$. Choose a Lagrangian $L_*$ of
$(\R^n\oplus{\R^n}^*,\omega)$ which is complementary to both
$\ell(t_0)$ and $L_0=\{0\}\oplus{\R^n}^*$; such Lagrangian always exists
(see for instance \cite[Corollary~3.2.9]{MPT}). Consider
the extended index form $I^\#_t$ defined in \eqref{eq:sustI}
corresponding to the choice of the bilinear form $\Theta=\phi_{L_1,L_0}(L_*)$,
where $L_1=\R^n\oplus\{0\}$. With such a choice, we have that
$I^\#_t$ is nondegenerate on $\mathcal H_t^\#$ for $t$ near $t_0$
(see formula \eqref{eq:KerIsust}).

Using Proposition~\ref{thm:indKsustindK}, for $t\ne t_0$ sufficiently
close to $t_0$ we have:
\begin{equation}\label{eq:fim1}
i(t)=n_-\Big(I^\#_t\big\vert_{\Kt_t^\#}\Big)-n_-\Big(\phi_{L_1,L_0}(L_*)
-\Phi_{L_1,L_0}\big(\ell(t)\big)\Big).
\end{equation}

Using Corollary~\ref{thm:corelementary},
Proposition~\ref{thm:IdComp}   and
Corollary~\ref{thm:cornegdef}, the jump of the
function $n_-\left(I^\#_t\big\vert_{\Kt_t^\#}\right)$ as $t$ passes through
$t_0$ equals the dimension of $\Kt_{t_0}\cap\St_{t_0}$, which by
 Corollary~\ref{thm:multfocinstred} is equal to the multiplicity of $t_0$ as a focal
instant for the reduced symplectic system. The sum of these multiplicities
as $t_0$ varies in $]a,b[$ equals $n_+(I\vert_{\St})$ by Corollary~\ref{thm:indIS}.

By Proposition~\ref{thm:calcMaslov2}, the jump of the function
$n_-\big(\phi_{L_1,L_0}(L_*) -\Phi_{L_1,L_0}(\ell(t))\big)$ as $t$ passes through $t_0$
is equal to 
$-\mu_{L_0}\big(\ell\vert_{[t_0-\varepsilon,t_0+\varepsilon]}\big)$ for $\varepsilon>0$
sufficiently small. Since $\mu_{L_0}$ is additive by concatenation,
the sum of these jumps equals minus the Maslov index of the Morse--Sturm
system. 

This concludes the proof for the case of a Morse--Sturm system
\eqref{eq:MS} whose focal instants
are isolated and such that $t=b$ is not focal for the reduced symplectic system
\eqref{eq:symplassoc}.

Consider now the more general case of a Morse--Sturm system 
for which $t=b$ is not focal for the associated reduced symplectic system.
Let $R_n:[a,b]\to\mathcal L(\R^n)$ be a sequence of real analytic
curves of $g$-symmetric linear endomorphisms of $\R^n$ that converges
uniformly to $R$ on $[a,b]$. Let $I^{(n)}$ be the index
form of the corresponding Morse--Sturm problem and denote by $\Kt^{(n)}$
the associated space defined as in \eqref{eq:defKS}. Then, $I^{(n)}$
converges to $I$ in the operator norm topology. Since  $I\vert_\St$ is nondegenerate
(see formula \eqref{eq:kerIX} and Proposition~\ref{thm:restrI})
and it is represented by a compact perturbation of
a negative isomorphism of $\St$ (Corollary~\ref{thm:Idefneg}),
it follows that, for $n$ sufficiently large,  $n_+(I^{(n)}\vert_\St)=n_+(I\vert_\St)$
(Remark~\ref{thm:remcasoparticular}).
Moreover, since $I\vert_\Kt$ is nondegenerate, $I\vert_\Kt$ is represented by
a compact perturbation of a positive isomorphism of $\Kt$ and
$\Kt^{(n)}$ converges%
\footnote{here we use the fact that $t=b$ is not focal for the reduced
symplectic system, as well as Corollary~\ref{thm:qFsurjective}.} 
to $\Kt$, by Remark~\ref{thm:remcasoparticular}  
we have $n_-(I\vert_\Kt)=n_-(I^{(n)}\vert_{\Kt^{(n)}})$ for $n$
sufficiently large.

The conclusion in the case that $t=b$ is not focal for the reduced
symplectic system follows from the
stability of the Maslov index by uniformly small perturbations of the coefficient
$R$ in \eqref{eq:MS} (see \cite[Theorem~5.2.1]{MPT}).

In the general case that $t=b$ may be focal for the reduced symplectic system,
the conclusion follows from the fact that the functions $i(t)$ and 
$n_+(I_t\vert_{\St_t})$ are {\em left-continuous\/} at $t=b$. 
The left-continuity of $n_+(I_t\vert_{\St_t})$ follows from Corollary~\ref{thm:indIS};
 the left-continuity of $i(t)$ follows from  Corollary~\ref{thm:cornegdef}
and from formula \eqref{eq:fim1}.
\end{proof}
\end{section}

\begin{section}{The Index Theorem for Symplectic Differential
Systems}\label{sec:indexthsympl}
In Subsection~\ref{sub:symplectic} we have defined the notion
of symplectic differential system, and we have seen that every such
system is isomorphic to a Morse--Sturm system (Proposition~\ref{thm:isoMSsympl}).
Moreover, we have seen that the notions of focal instants, multiplicity,
signature and index form are invariant by isomorphisms (Proposition~\ref{thm:isosympl}).
This suggests that it is possible to give a general version of the index theorem 
for symplectic differential systems; the purpose of  this section 
is to give the main definitions and to state the generalized index
theorem for symplectic systems with initial conditions.
We will use most of the notations introduced in Subsection~\ref{sub:symplectic};
the details of many of the results presented in this section
may be found in \cite[Section~2]{PT2}.

We consider a symplectic differential system in $\R^n$
of the form \eqref{eq:sympl},
and we consider the initial conditions:
\begin{equation}\label{eq:symplIC}
v(a)\in P,\quad \alpha(a)\vert_P+S\big(v(a)\big)=0,
\end{equation}
where $P\subset \R^n$ is a subspace and $S$ is a symmetric
bilinear form on $P$, considered as a map from $P$ to $P^*$.
The set $\ell_0\subset\R^n\oplus{\R^n}^*$ defined by:
\begin{equation}\label{eq:ell0}
\ell_0=\Big\{(v,\alpha):v\in P,\ \alpha\vert_P+S(v)=0\Big\}
\end{equation}
is a Lagrangian subspace of $(\R^n\oplus{\R^n}^*,\omega)$;
conversely, every Lagrangian subspace $\ell_0$
of $(\R^n\oplus{\R^n}^*,\omega)$
defines uniquely a subspace $P\subset\R^n$ and a symmetric bilinear
form on $P$ such that \eqref{eq:ell0} holds. We will
say that $v$ is an {\em $(X,\ell_0)$-solution\/} if $v$ is an $X$-solution
such that $(v(a),\alpha_v(a))\in\ell_0$. In analogy with \eqref{eq:spaces},
we now define:
\begin{equation}\label{eq:novoV}
\mathbb V=\Big\{v: v\ \text{is an $(X,\ell_0)$-solution}\Big\}.
\end{equation}
The notions of focal instant, multiplicity, signature and focal index
for the pair $(X,\ell_0)$ are given  in Definition~\ref{thm:deffocal},
where the space $\mathbb V$ is now redefined in \eqref{eq:novoV}.

As in the case of  semi-Riemannian geodesics, we need the following
nondegeneracy assumption on the initial conditions for
the symplectic differential system:
\begin{defin}\label{thm:defICsympl}
A pair $(X,\ell_0)$ where $X$ is the coefficient matrix of
a symplectic differential system and $\ell_0$ is a Lagrangian
subspace of $(\R^n\oplus{\R^n}^*,\omega)$ is said to be a {\em set of
data for the symplectic differential problem\/} if the symmetric
bilinear form $B(a)^{-1}$ in $\R^n$ is nondegenerate on the subspace
$P$ associated to $\ell_0$.
\end{defin}

Let $(X,\ell_0)$ be a set of data for the symplectic differential
problem; 
for each $t\in[a,b]$, the subspace $\ell(t)\subset\R^n\oplus{\R^n}^*$ given
by:
\[\ell(t)=\Big\{\big(v(t),\alpha_v(t)\big): v\in\mathbb V\Big\}\]
is Lagrangian. So, we get a $C^1$-curve $\ell$ in the Lagrangian Grassmannian
$\Lambda$; recalling the notations of Subsection~\ref{sub:Maslovcurve} and
setting $L_0=\{0\}\oplus{\R^n}^*$, it is
easily seen that $\ell(t)\in\Lge(L_0)$ if and only if $t$ is a focal instant.

Given the nondegeneracy assumption in Definition~\ref{thm:defICsympl},
it is possible to see that there exists $\varepsilon>0$ such that there
are no focal instants in $\left]a,a+\varepsilon\right]$. We can therefore give the
following definition:
\begin{defin}\label{thm:defMaslovsympl}
If $t=b$ is not a focal instant, the {\em Maslov index\/} $\imaslov(X,\ell_0)$
of the pair $(X,\ell_0)$ is defined as:
\begin{equation}
\imaslov(X,\ell_0)=\mu_{L_0}\big(\ell\vert_{[a+\varepsilon,b]}\big),
\end{equation}
where $\varepsilon>0$ is chosen in such a way that there are no
focal instants in $]a,a+\varepsilon]$.
\end{defin}

Proposition~\ref{thm:imaslovifoc} can be generalized to symplectic
systems.

The {\em index form\/} $I_{(X,\ell_0)}$ associated to the symplectic differential
problem is the bounded symmetric bilinear form on the Hilbert space $\mathcal H$ given
in \eqref{eq:defH} defined by:
\begin{equation}\label{eq:defIXell0}
I_{(X,\ell_0)}(v,w)=\int_a^b\Big[B(\alpha_v,\alpha_w)+C(v,w)\Big]\;\mathrm dt
-S\big(v(a),w(a)\big).
\end{equation}
Recalling Definition~\ref{thm:defequiv}, we now give the following 
definition of isomorphisms for symplectic differential systems
with initial data:
\begin{defin}\label{thm:defisoinitial}
The pairs $(X,\ell_0)$ and $(\tilde X,\tilde\ell_0)$ of data for the symplectic
differential problem are said to be {\em isomorphic\/} if there exists
an isomorphism $\phi_0$ between $X$ and $\tilde X$ such that
$\phi_0(a)(\ell_0)=\tilde\ell_0$.
\end{defin}

\noindent
Proposition~\ref{thm:isosympl} generalizes {\em mutatis mutandis\/}
to the case of isomorphisms of pairs $(X,\ell_0)$; moreover, isomorphic
pairs have the same Maslov index (see \cite[Proposition~2.10.2]{PT2}).

Using Proposition~\ref{thm:isoMSsympl}, we have the following index theorem
for symplectic systems:
\begin{teo}\label{thm:indextheoremsympl}
Let $(X,\ell_0)$ be a smooth set of data for the symplectic differential
problem in $\R^n$, with $k=n_-(B)$. Let $Y_1,\ldots,Y_k:[a,b]\to\R^n$
be  smooth maps such that, for each $t\in[a,b]$, $Y_1(t),\ldots,Y_k(t)$
for a basis of a subspace $D_t\subset\R^n$ on which $B(t)^{-1}$ is
negative definite. Consider the following two closed subspaces of $\mathcal H$
(see \eqref{eq:defH}):
\begin{equation}\label{eq:defKSsympl}
\begin{split}
\Kt=\Big\{v\in\mathcal H:&\;\alpha_v(Y_i)\in H^1([a,b];\R)\
\text{and}
\\&
\alpha_v(Y_i)' =B(\alpha_v,\alpha_{Y_i})+C(v,Y_i),\ \ \forall\,i=1,\ldots,k\Big\};\\
\St=\Big\{v\in H^1_0&\;([a,b];\R^n):v(t)\in D_t,\ \forall \, t\in [a,b]\Big\}.
\end{split}
\end{equation}
Then, if $t=b$ is not focal, we have:
\begin{equation}\label{eq:indexsympl}
\imaslov(X,\ell_0)=n_-\left(I_{(X,\ell_0)}\big\vert_\Kt\right)-
n_+\left(I_{(X,\ell_0)}\big\vert_\St\right)-n_-\left(B(a)^{-1}\big\vert_P\right),
\end{equation}
where all the terms in the above equality are finite integer numbers.
\end{teo}
\begin{proof}
It follows directly
from Proposition~\ref{thm:isoMSsympl} and the proof of Theorem~\ref{thm:indexth}.
\end{proof}
It is easy to see that the space $\Kt$ depends only on the family
of subspaces $\{D_t\}_{t\in[a,b]}$, and not on the particular choice of
a basis $Y_1,\ldots,Y_k$. Moreover, the spaces $\Kt$ and $\St$ are
$I_{(X,\ell_0)}$-orthogonal.

Also in this context it is possible to determine a {\em reduced
symplectic system\/} associated to the choice  of
the vector fields $Y_i$. The formula of this reduced system is the same as
\eqref{eq:symplassoc}, where the matrices $\mathcal B$, $\mathcal C$ and
$\mathcal I$ are now given by:
\begin{equation}\label{eq:calobjsympl}
\mathcal B_{ij}=B^{-1}(Y_i,Y_j),\quad \mathcal C_{ij}=\alpha_{Y_j}(Y_i),
\quad
\mathcal I_{ij}=B(\alpha_{Y_i},\alpha_{Y_j})+C(Y_i,Y_j).
\end{equation}
Observe that the reduced symplectic system is always considered
with initial conditions $f(a)=0$ regardless of the initial conditions
considered for the original symplectic system.

Many of the results of Subsection~\ref{sub:reduced} (Lemma~\ref{thm:systred},
Corollary~\ref{thm:multfocinstred}, Proposition~\ref{thm:restrI},
Corollaries~\ref{thm:indIS} and \ref{thm:Idefneg}) generalize to
this context. In particular, if $t=b$ is not focal for the reduced
symplectic system, then $\mathcal H=\Kt\oplus\St $ and the
term $n_+\left(I_{(X,\ell_0)}\vert_{\St}\right)$ in formula
\eqref{eq:indexsympl} can be computed as the sum of the multiplicities
of the focal instants of the reduced symplectic system in $]a,b[$.

\begin{rem}\label{thm:remregularity}
Observe that the proof of Theorem~\ref{thm:indextheoremsympl} is valid
under a weaker assumption on the regularity of the coefficients
of the symplectic system and of the fields $Y_i$. More precisely, if
$t=b$ is not focal for the reduced symplectic system, our proof works
in the case that $A$ is of class $C^1$, $B$ is of class $C^2$, $C$ is continuous
and the $Y_i$'s are of class $C^2$. In the general case, one has to assume that
$A$ is of class $C^2$, $B$ is of class $C^3$, $C$ is of class $C^1$ and
the $Y_i$'s are of class $C^3$. It is known to the authors that 
a direct proof of Theorem~\ref{thm:indextheoremsympl} (that does not
use Proposition~\ref{thm:isoMSsympl}), technically more involved than
the one presented in this paper, shows that the regularity assumption
can be weakened even further. Namely, if $t=b$ is not
focal for the reduced symplectic system, Theorem~\ref{thm:indextheoremsympl} is valid
under the assumption that $A$ and $B$ are of class $C^1$, $C$ is continuous and
the $Y_i$'s are of class $C^2$; in the general case one needs the assumption that
also $C$ is of class $C^1$. 
\end{rem}
\end{section}



\end{document}